\documentclass[12pt]{article}
\usepackage{amsmath}
\usepackage{amsthm}
\usepackage{subfigure}
\usepackage[xetex]{graphicx}
\usepackage{epstopdf}
\usepackage{amssymb,enumerate}
\usepackage{graphicx}
\usepackage{latexsym}
\usepackage{xcolor}

\newtheorem{thm}{Theorem}[section]
\newtheorem{cor}{Corollary}[section]
\newtheorem{lem}{Lemma}[section]
\newtheorem{prop}{Proposition}[section]

\newtheorem{rem}{Remark}[section]

\newcommand{\be}{\begin{equation}}
\newcommand{\ee}{\end{equation}}

\newcommand{\norm}[1]{\left\Vert#1\right\Vert}

\makeatletter

\newcommand{\Rmnum}[1]
{expandafter\@slowromancap\romannumeral#1@}
\makeatother
\numberwithin{equation}{section} 
\numberwithin{figure}{section} 
\allowdisplaybreaks[4]

\begin{document}
\title{Sharp Thresholds of Instability   for the  Biharmonic Nonlinear Schr\"odinger Equations}
\author{Shanshan Fu$^{1,2}$,\quad Shihui Zhu$^1$\thanks{Corresponding author's \textit{E-mail address}: shihuizhumath@163.com; shihuizhumath@sicnu.edu.cn.  The 16-digit ORCID of the corresponding author: 0000-0001-7493-9123.}\\
{\small \textit{$^1$School of Mathematical Sciences, Sichuan Normal University,}}
\\ {\small \textit{ Chengdu 610066, PR China}}\\
{\small \textit{$^2$ School of Science, Xihua University,Chengdu 610039, PR China}}}

\date{}
\maketitle
{\begin{center}
\begin{minipage}{6.0in}
{\textbf{Abstract:}   We consider the time-dependent biharmonic nonlinear Schr\"odinger equations with and without a mixed dispersion. We find the sharp threshold of blow-up and global existence for the biharmonic nonlinear Schr\"odinger equations without a mixed dispersion by establishing the modified cross-constraint variational method. Building on this work, we prove that the standing waves are strongly unstable without further assumptions on the frequency.  The main new ingredients of which require us to provide a comparison of two variational problems, improved the method in Zhang [Comm PDEs, 2005, 30: 1429-1443]. 
  In the latter part of the paper, we study the biharmonic nonlinear Schr\"odinger equations with a mixed dispersion, which destroys the scaling invariance. The single variational approach is utilized to obtain the sharp threshold of blow-up and global existence, and strong instability of standing waves. }\\

{\textbf{Key words:} fourth-order nonlinear Schr\"odinger equation; mixed dispersion;  modified cross-constraint variational method; global existence; strong instability.}\\
 {\textbf{2020 Mathematics Subject Classification:}  35B35\ \ 35Q55 \ \ 35J50.}  
\end{minipage}
\end{center}

\vspace{0.3cm}
 \section{Introduction}
 \indent
 
We consider the Cauchy problem of the biharmonic nonlinear Schr\"odinger equations
\begin{equation}\label{eq1} 
  \begin{cases}
  &i\phi_t - \Delta^{2} \phi + \mu \Delta\phi + \sigma  |\phi|^{p-1}\phi =0, \quad \text{in}  \,\ \mathbb{R}^{N}\times  \mathbb{R^+},\\
  &\phi(x, 0)=\phi_{0}(x), \quad \text{in}  \,\ \mathbb{R}^{N},
  \end{cases}
\end{equation}
where $\phi=\phi(x,t) : \mathbb{R}^N\times\mathbb{R}^+\to\mathbb{C}$ is a complex function, $\Delta^2 = \Delta (\Delta)$ is the biharmonic operator, and $\Delta=\sum_{j=1}^{N}\frac{\partial^2}{\partial x_j^2}$ is the Laplace operator. The parameter $\mu \geqslant 0$ represents the strength of the second-order dispersive term, and $\sigma \in \mathbb{R}$ coefficient of  nonlinear interaction, where it is  said to be defocusing when $\sigma < 0 $ and focusing when $\sigma > 0$.

In this paper, we concentrate on the focusing case $\sigma > 0$, under the assumptions that the spatial dimension $ N \geqslant 2$, the dispersion coefficient $\mu \geqslant 0$, and the nonlinear exponent $p$ satisfies $\frac{N+8}{N}\leqslant  p<\frac{N+4}{ (N-4)^{+} }$, where $\frac{N+4}{ (N-4)^{+} }=+\infty$ for $N \leqslant 4$ and $\frac{N+4}{ (N-4)^{+} }=\frac{N+4}{N-4}$ for $N > 4$, which reflects the critical Sobolev embedding threshold.  Ben-Artzi et al \cite{BenArtzi2000},  Kenig et al \cite{kenig1991}, Pausader \cite{Pausader2007}, and Cui and Guo \cite{CuiGuo2007}, etc., have studied the well-posedness theory for the Cauchy problem (\ref{eq1}) in various work spaces (see the book \cite{Fibich2015} and references therein). In particular, if the initial data
 $\phi_{0}(x) \in  H^2(\mathbb{R}^N)$, then there exists a unique solution $\phi(x,t)\in C([0,T);H^2(\mathbb{R}^N))$, and $[0, T)$ is the maximal existence interval. 
 Furthermore, for all $t\in [0, T)$,  $\phi(x,t)$ enjoys two fundamental conservation laws.
\begin{itemize}
  \item[(i)] Mass
\begin{equation}\label{eqm}
  M(\phi(t))=\int|\phi(t, x)|^2dx= M(\phi_0).
\end{equation}
  \item[(ii)] Energy
  \begin{equation}\label{eqe1}
 E(\phi(t))=\frac{1}{2}\int \left|\Delta \phi(t, x)\right|^2 + \mu |\nabla \phi(t, x)|^2dx-\frac{\sigma}{p+1}\int \left|\phi(t, x)\right|^{p+1}dx= E(\phi_0).
\end{equation}
\end{itemize} 
 
 Essentially, Eq.~\eqref{eq1} is a class of fourth-order nonlinear Schr\"odinger equations (FONLSs), proposed by 
 Karpman and Shagalov \cite{KS2000}, considering the role of small fourth-order dispersion. Now, the FONLSs is becoming one of the most important models in nonlinear optics, and it can model the evolution of intense laser beams in media where the standard paraxial approximation breaks down, and higher-order dispersion needs to be considered, and  the description of ultra-short pulse propagation in optical fibers, where the fourth-order term accounts for corrections to the standard group velocity dispersion (see \cite{Fibich2002,Fibich2015} and references therein).  
 Actually, there are two crucial  questions in nonlinear Schr\"odinger equation theory.
 
\begin{itemize}
\item Do solutions undergo singularity formation (such as blow-up or collapse) within a finite time, or do they scatter?
\item Are there finite-energy waveguide solutions of the form $\phi(x,t)=e^{i\omega t} Q(x)$, and if they exist, what is their stability property?
\end{itemize}
As anticipated, these two problems have been the focus of significant research in mathematics, and many ingenious methods and various outstanding results are emerging, including the virial identity and the general sufficient conditions for the existence of blow-up solutions \cite{HolmerR2010, Weinstein1983}, cross-constraint variational method and sharp threshold for blow-up and global existence \cite{Zhang2002}, $I$-method and scattering (non-existence of blow-up) \cite{CKSTT2002},  the spectrum method and the existence of small super-critical mass blow-up solutions without $|x|\phi_0\in L^2$  \cite{Merle2005}, the concentration compactness regality method and sharp threshold of scattering versus blow-up \cite{KenigMerle2006IM}, concentration compactness principle and orbital stability \cite{CL1982}, linear and nonlinear stability frame \cite{GSS1987}, \cite{Lin2022} etc. 

However, most of the above methods fail for the FONLSs. As might be expected, the  
FONLSs have received considerable mathematical attention, including the ground and bounded solutions \cite{FJMM2022, LW2024, Ma2023}, global existence \cite{CuiGuo2007, Fibich2002, HN2025, ZhuZhangYang2011}, scattering \cite{MZ2016, Pausader2007b, PX2013, SS2020, YYZ2023, YYZ2024}, existence of blow-up solutions \cite{Boulenger2017, Fibich2003, K2022}, stability of standing waves \cite{BCGJ2019TAMS, BCGJ2019, LuoZhengZhu2023, NP2015} etc. Although there are many excellent results for the FONLSs, to study the existence of blow-up solutions for Eq.~\eqref{eq1} is of particular interest and challenge. For one thing, on the existence of blow-up solutions, Fibich \cite{Fibich2002} performed numerical simulations for the fourth-order equation and observed self-similar blow-up profiles in the mass-supercritical regime.  The authors of the paper \cite{ZhuZhangYang2011} rigorously studied the structure of singularities by assuming the existence of blow-up solutions and analyzing their concentration behavior. Until 2017, 
a breakthrough was made by Boulenger and Lenzmann \cite{Boulenger2017}, who developed a generalized virial-type identity tailored for fourth-order equations, leading to a differential inequality controlling the time evolution of a localized virial quantity $\mathcal{M}_{\varphi_R}[u(t)]$. This inequality incorporates higher-order derivative terms and spatial cutoff functions, and plays a crucial role in proving finite-time blow-up under suitable energy and mass conditions. But how to distinguish the initial data to blow-up regime and global regime is still unsolved, that is, a sharp threshold of blow-up and global existence for Eq.~\eqref{eq1}.

\newpage
The main purpose of this paper is to prove the following result.
\begin{thm}[Sharp threshold I: $\mu = 0$]\label{th1}
Let $\mu = 0$ and $N \geqslant  2$. Assume initial data $\phi_0(x)\in H^2$,
$$E(\phi_0)+\frac{\omega}{2}M(\phi_0)< D_{\mathcal{M}^0_0}, $$
where $D_{\mathcal{M}^0_0}$ is the minimum of a variational problem defined in Section 3. Then, we have the following sharp threshold of blow-up and global existence for the Cauchy problem \eqref{eq1} in the following sense. 
\begin{itemize}
\item [(i)] If $\phi_0(x)$ is radial, $ \frac{N+8}{N}< p < \min\left\{\frac{N+4}{(N-4)^{+}}, 9\right\} $ and satisfies
$$  \int|\Delta \phi_0|^{2}dx+\omega\int|\phi_0|^{2}dx<\sigma\int|\phi_0|^{p+1}dx  \  {\rm and }\  \int|\Delta \phi_0|^{2}dx < \frac{N\sigma(p-1)}{4(p+1)}\int|\phi_0|^{p+1}dx.$$
Accordingly, the solution $\phi ( t, x)$ to the Cauchy problem \eqref{eq1} will exhibit blow-up in a finite time.
\item [(ii)] If $ \frac{N+8}{N}< p < \frac{N+4}{(N-4)^{+}} $ such that  $$ \int|\Delta \phi_0|^{2}dx+\omega\int|\phi_0|^{2}dx<\sigma\int|\phi_0|^{p+1}dx  \  {\rm and }\  \int|\Delta \phi_0|^{2}dx > \frac{N\sigma(p-1)}{4(p+1)}\int|\phi_0|^{p+1}dx,$$
or $$\int|\Delta \phi_0|^{2}dx +\omega\int|\phi_0|^{2}dx>\sigma\int|\phi_0|^{p+1}dx.$$
\end{itemize}
Then, the corresponding solution $\phi ( t, x)$  to Cauchy problem \eqref{eq1} exists globally for all $t$ in the interval $[0, \infty)$.
\end{thm} 
 \begin{rem}\label{remhist}
Due to that the classical virial argument does not apply to $\Delta^{2}$, the rigorous proof of the existence of  finite-time blow-up  solutions for the FONLS \eqref{eq1} 
is an open problem by Baruch, Fibich and Mandelbaum
in \cite{BFM2010}, where  the existence of
blow-up solutions was supported only by the numerical studies.
Until 2017, Boulenger and Lenzmann in \cite{Boulenger2017} settled this question by  defining a new 
localized virial estimate under the constraint of radial initial data with negative energy. When
$\mu=0$ and $p>\frac{N+8}{N}$, they also allow data below the ground state threshold.
 Dinh in \cite{Dinh2019,Dinh2021} further proved a sharp dichotomy between scattering and blow-up below that
threshold for radial data for $\mu>0$.
\newpage
However, by improving the cross-constraint variational arguments proposed by Zhang in \cite{Zhang2005}, we obtain the sharp threshold of blow-up and global existence for  the FONLS \eqref{eq1} in Theorem \ref{th1}.
Here, the global existence aspect, we do not need the radial assumption on initial data, and the initial energy may be positive. More of importance, we can use this work frame to 
give the strong instability of Theorem \ref{th5.1}, with no constraint on $\omega$.
\end{rem}

Here, it is worth noting that the global existence conclusions established in Theorem \ref{th1} are extendable to the critical case ($p=\frac{N+8}{N}$), leading to the following corollary.
\begin{cor}[Global existence: $\mu = 0$ and $p=\frac{N+8}{N}$]\label{cor1}
Let $N \geqslant  2$, $\mu = 0$ and $p=\frac{N+8}{N}$. Assume  that  $\phi_0(x)\in H^2$ is such that 
$$E(\phi_0)+\frac{\omega}{2}M(\phi_0)< D_{\mathcal{M}^0_0},$$
where $D_{\mathcal{M}^0_0}$ is the minimum of a variational problem defined in Section 3, and satisfies one of the following conditions
\begin{itemize}
\item [(i)] $$\int|\Delta \phi_0|^{2}dx+\omega\int|\phi_0|^{2}dx<\sigma\int|\phi_0|^{p+1}dx\,\
\text{and}\,\
8\int|\Delta \phi_0|^{2}dx > \frac{2N\sigma(p-1)}{p+1}\int|\phi_0|^{p+1}dx.$$
\item [(ii)] $$\int|\Delta \phi_0|^{2}dx +\omega\int|\phi_0|^{2}dx>\sigma\int|\phi_0|^{p+1}dx.$$
\end{itemize}
Then, the solution $\phi ( t, x)$ corresponding to the Cauchy problem \eqref{eq1} exists globally for all $t$ in the interval $[0, \infty)$.
\end{cor}

Next, we address the problem of identifying the sharp threshold between blow-up and global existence for the Cauchy problem \eqref{eq1} in the case $\mu>0$. In this regime, the lack of scaling symmetry becomes a principal challenge, making the cross-constraint variational method applied in Theorem \ref{th1} unusable.
 We establish a new variational method based on the profile decomposition theory. This method not only provides the sharp threshold for the blow-up and global existence to Eq.~\eqref{eq1}, but also can be applied to prove the strong instability of standing waves. The main results are as follows.
 \newpage
 \begin{thm}[Sharp threshold II: $\mu > 0$]\label{thmdi1}
        Let $\mu > 0$, $N \geqslant  2$  and $\phi_0\in H^2$ satisfies
$$E(\phi_0)+\frac{\omega}{2}M(\phi_0)< D_{\mathcal{M}^{\mu}_0},$$
where $D_{\mathcal{M}^\mu_0}$ is the minimum of a variational problem defined in Section 3. 
Then, the following conclusions hold.
\begin{itemize}
\item [(i)]  For $ \frac{N+8}{N}\leqslant p < \min\left\{\frac{N+4}{(N-4)^{+}}, 9\right\} $, if  the initial data $\phi_0(x)$ is radial and satisfies
$$8\int|\Delta \phi_0|^{2}dx+4\mu\int|\nabla \phi_0|^{2}dx < \frac{2N\sigma(p-1)}{p+1}\int|\phi_0|^{p+1}dx,$$
then, the solution $\phi(t,x)$ associated with the Cauchy problem \eqref{eq1} blows up in a finite time.
\item [(ii)]  For $\frac{N+8}{N} < p<\frac{N+4}{(N-4)^{+}}$, if  the initial data is such that 
$$8\int|\Delta \phi_0|^{2}dx+4\mu\int|\nabla \phi_0|^{2}dx > \frac{2N\sigma(p-1)}{p+1}\int|\phi_0|^{p+1}dx,$$
then the solution $\phi(t,x)$ associated with the Cauchy problem \eqref{eq1} exists for all time.
\end{itemize}
\end{thm}


Motivated by the variational and invariant manifold framework of Zhang \cite{Zhang2005} and building upon the virial-type blow-up method of Boulenger-Lenzmann \cite{Boulenger2017}, in this work, we aim to fill this gap by establishing sharp threshold conditions for blow-up and global existence, and proving the strong instability of standing waves. Specifically, by developing variational approaches integrating localized virial identity to analyze \eqref{eq1}, we construct constrained variational problems and invariant manifolds along the evolution flow, enabling us to characterize the sharp threshold conditions for global existence and finite-time blow-up. By incorporating a refined virial identity inspired by \cite{Boulenger2017}, we establish blow-up results in the absence of classical virial techniques. 
 
On the other hand,  Eq.~\eqref{eq1}  admits a special class of solutions known as standing waves, which take the form $\phi(x,t) = e^{i\omega t} u(x)$ with frequency $\omega > 0$. Substituting this ansatz into equation \eqref{eq1} leads to the associated elliptic equation:
\begin{equation}\label{eq64}
\Delta^2 u + \omega u - \mu\Delta u- \sigma |u|^{p-1} u = 0,
\end{equation}
which has an infinite number of solutions in $H^2$. Among these, there is exactly one solution of
minimal energy, called the ground state solution $Q$. The general solution of the above equation is called the
bound state solution,  plays a fundamental role in understanding the long-time behavior of FONLS and serves as a natural candidate for threshold analysis between global existence and blow-up. In this paper, we further investigate the strong instability of standing wave solutions to the FONLS equation.

  \begin{thm}[Strong instability: $\mu \geqslant 0$]\label{th5.1}
 For $\mu \geqslant 0$ and $N \geqslant  2$. Let $p$ be odd integer satisfying $ \frac{N+8}{N}\leqslant p < \min\left\{\frac{N+4}{(N-4)^{+}}, 9\right\} $. Then the standing waves of Eq.~\eqref{eq1} are strongly unstable. That is, for any $\varepsilon >0$ and the radial  minimizer $u$ of \eqref{eqN}, there exists $\phi_0\in H^2$ with $\|\phi_0-u\|_{H^2}<\varepsilon $ such that 
the corresponding solution $\phi (t, x)$ of Eq.~\eqref{eq1} blows up in a finite time.
\end{thm}
\begin{rem}
Although the standing wave solutions are strongly unstable for both $\mu = 0$ and $\mu > 0$, it is necessary to treat the two cases separately. This is due to essential differences in the equation structure, the attainability of the variational problem \eqref{eqN}, and the characterization of sharp blow-up thresholds.
\end{rem}
\begin{rem}
  It is worth emphasizing that the analysis of strong instability for the standing wave solutions to equation \eqref{eq1} is primarily based on the method of Zhang \cite{Zhang2005}. However, we remove a key condition in their work, namely the assumption that the frequency $\omega$ satisfies $d_{\mathcal{M}} \geqslant d_{\mathcal{N}}$, where $d_{\mathcal{M}}$ and $d_{\mathcal{N}}$ are the minimal values associated with two constrained variational problems. In our setting, we denote these quantities by $D_{\mathcal{M}_{-}^{0}}$ and $D_{\mathcal{M}^{0}_{0}}$, respectively, and establish the ordering $D_{\mathcal{M}^{0}_{0}} \leqslant D_{\mathcal{M}_{-}^{0}}$ by rigorous proof in Theorem \ref{lem12}, thereby justifying and improving upon the condition used in \cite{Zhang2005}.
\end{rem}
\begin{rem}
 We should point out that although the strong instability  of standing waves for Eq.~\eqref{eq1} has been investigated by  Bonheure et al \cite{BCGJ2019}, we prove that the standing waves for Eq.~\eqref{eq1} 
are strongly unstable  under $\mu=0$, and $\mu>0$, by developing distinct argumentation. More importantly, we prove the constrained minimization problem (\ref{eqN}) can be attained, but this is not proved in Proposition 2.1,  which is crucial to proving the strong instability in \cite{BCGJ2019}.
\end{rem}
At the end of this section, we provide an overview of the structure of the paper. Section 2 presents preliminary results. Section 3 constructs the invariant manifolds. Section 4 establishes the sharp threshold for global existence and blow-up. Section 5 proves the strong instability of standing waves.

\section{Notations and Preliminaries}
In this section, we review several known results and establish some preliminary facts that will be essential in the subsequent analysis.
 For convenience, we denote $L^p(\mathbb{R}^N)$, $\|\cdot\|_{L^p(\mathbb{R}^N)}$, $H^2(\mathbb{R}^N)$, and $\int_{\mathbb{R}^N}\cdot dx$ by $L^p$, $\|\cdot\|_{p}$, $H^2$, and $\int \cdot dx$, respectively.
Furthermore, we shall write $X \lesssim  Y$ to denote that $X \leqslant CY$ holds with constant $C>0$.

The work space for Eq.~\eqref{eq1} is $H^2=\{u\in L^2| \nabla u\in L^2 \,\ \text{and} \,\ \Delta u\in L^2\}$, equipped the \[\|u\|_{H^2}^2=\|u\|_2^2+\|\nabla u\|_2^2+\|\Delta u\|_2^2.\] One can check that 
\begin{equation}\label{eq76}
C_1(\|u\|_2^2+\|\Delta u\|_2^2)\leqslant \|u\|_{H^2}^2\leqslant  C_2(\|u\|_2^2+\|\Delta u\|_2^2),
\end{equation}
where $C_1$ and $C_2$ are two positive constants. 
As a starting point, we recall the local well-posedness for the Cauchy problem \eqref{eq1} in $H^2$ (see \cite{Pausader2007} and the references).
\begin{lem}[Local Well-posedness \cite{Pausader2007}]\label{prop2.1}
Given any initial data \( \phi_0 \in H^2 \), any \( p \in (1, \frac{N+4}{N-4}) \) when \( N \geqslant 5 \), and any \( p > 1 \) when \( N \leqslant 4 \), there exists \( T > 0 \) and a unique solution \( \phi \in C(H^2, [0, T]) \) of \eqref{eq1}. 
And the solution preserves the quantities of mass and energy, in the sense that
\begin{equation}\label{4.2}
M(\phi(t)) = M(\phi_0) \quad \text{and} \quad E(\phi(t)) = E(\phi_0) \quad \text{for all } t \in [0,T],
\end{equation}
with $M$ and $E$ given by \eqref{eqm} and \eqref{eqe1}.
In addition, if $T^*<+\infty$ denotes the maximal time of existence, then
\begin{equation*}
\lim_{t\to T^*} \|\phi(t)\|_{H^2} = +\infty.
\end{equation*}
\end{lem}

In preparation for the analysis of finite-time blow-up and standing wave instability for Eq.~\eqref{eq1}, 
we recall the localized virial inequality due to Boulenger and Lenzmann~\cite{Boulenger2017}. Specifically,
let $\varphi : \mathbb{R}^N \to \mathbb{R}$ be radial, with $\nabla^j \varphi \in L^\infty(\mathbb{R}^N)$ for $1 \leqslant j \leqslant 6$, 
and satisfying
\[
\varphi(r) =
\begin{cases}
r^2/2, & r \leqslant 1, \\
\text{const}, & r \geqslant 10,
\end{cases}
\]
as well as $\varphi''(r) \leqslant 1$ for all $r \geqslant 0$.
Given $R > 0$, we rescale by defining $\varphi_R(r) := R^2 \varphi(\frac{r}{R})$.
Then the localized virial of $u\in H^2$ is defined by
\begin{equation}\label{eq75}
    \mathcal{M}_{\varphi_R}[u]:=\langle u,-i(\nabla\varphi_R\cdot\nabla+\nabla\cdot\nabla\varphi_R)u\rangle = 2\,\mathrm{Im}\int \overline{u}\nabla\varphi_R\cdot\nabla u dx.
\end{equation}
\begin{lem}[Time Evolution of $\mathcal{M}_{\varphi_R}$ \cite{Boulenger2017}]\label{lem11}
 Let $N\geqslant 2$ and $R>0$. Suppose that $u\in C([0,T);H^2)$ is a radial solution of Eq.~\eqref{eq1}. Then, for any $t \in [0, T)$, 
\begin{align*}
\frac{d}{dt} \mathcal{M}_{\varphi_R}[u]\leqslant &2N(p-1) E[u_0] - (N(p-1) - 8)\|\Delta u(t)\|_{2}^2 - (N(p-1) - 4)\mu\|\nabla u(t)\|_{2}^2 + X_\mu[u(t)]\\
& + O\left(R^{-4} + R^{-2}\|\nabla u(t)\|_{2}^2 + R^{-\frac{(p-1)(N-1)}{2}}\|\nabla u(t)\|_{2}^{\frac{p-1}{2}} + |\mu|R^{-2}\right)\\
= & H(u)+O\left(R^{-4} + R^{-2}\|\nabla u(t)\|_{2}^2 + R^{-\frac{(p-1)(N-1)}{2}}\|\nabla u(t)\|_{2}^{\frac{p-1}{2}} + |\mu|R^{-2}\right),
\end{align*}
where
\begin{equation*}
  X_\mu[u] \leqslant
\begin{cases}
0, & \mu \geqslant 0,\\
K|\mu|\|\nabla u(t)\|_{2}^2, & \mu < 0,
\end{cases}
\end{equation*}
with constant $K > 0$ depending only on $\varphi$.
\end{lem}

Now we introduce the profile decomposition of bounded sequences in $H^2$, proposed by Zhu, Zhang, and Yang \cite{ZhuZhangYang2010}, which has been successfully applied to find new sharp Gagliardo-Nirenberg (G-N) inequality \cite{ZhuZhangYang2010}, and to prove orbital stability of standing waves for nonlinear Schrödinger equations (see \cite{ZhuJEE2017}).
\begin{lem}[Profile Decomposition \cite{ZhuZhangYang2010}]\label{lemdi51}
Let $\{u_n\}_{n = 1}^{+\infty}$ be a bounded sequence in $H^2$. Then there exists a subsequence of $\{u_n\}_{n = 1}^{+\infty}$ $($ for convenience, we still denote it as $\{u_n\}_{n = 1}^{+\infty}$ $)$, a family $\{x_n^j\}_{j = 1}^{+\infty}$ of sequences in $\mathbb{R}^N$ and a family $\{U^j\}_{j = 1}^{+\infty}$ of functions in $H^2$ such that
\begin{enumerate}[(i)]
    \item for every $k\neq j$
    \[
    |x_n^k - x_n^j| \to \infty \quad \text{as }n \to \infty;
    \]
    \item for every $l\geqslant1$ and every $x\in\mathbb{R}^N$, $u_n$ can be decomposed as
    \[
    u_n(x)=\sum_{j = 1}^{l}U^j(x - x_n^j)+u_n^l(x),
    \]
    with
    \[
    \lim_{l\to \infty}\limsup_{n\to \infty}\norm{u_n^l}_q = 0,\quad\text{for every }q\in\left(2,\frac{2N}{(N-4)^+}\right).
    \]
\end{enumerate}
Moreover,  the norms of $ \{u_n\}$ have the following decomposition
\begin{align}
    \norm{u_n}_2^2&=\sum_{j = 1}^{l}\norm{U^j}_2^2+\norm{u_n^l}_2^2+o_n(1),\label{eq:1.8}\\
    \norm{\Delta u_n}_2^2&=\sum_{j = 1}^{l}\norm{\Delta U^j}_2^2+\norm{\Delta u_n^l}_2^2+o_n(1),\label{eq:1.08}
\end{align}
for every $l \geqslant 1$, and $\lim\limits_{n\to\infty}o_n(1)=0$.
\end{lem}

\begin{prop}[Sharp G-N Inequality \cite{ZhuZhangYang2011}]\label{prop2.4}
Let $N \geqslant 1$ and  
define the variational problem
\begin{equation*}
\alpha := \inf\limits_{0 \neq v \in H^2} I^{N,p}(v) = 
\inf\limits_{0 \neq v \in H^2} \frac{\left( \int |v|^2 dx \right)^{\frac{\theta_1}{2}}\left( \int |\Delta v|^2 dx \right)^{\frac{\theta_2}{2}}}{\int |v|^{p + 1} dx},
\end{equation*}
where $\theta_1=p+1-\frac{(p-1)N}{4}$ and $\theta_2=\frac{(p-1)N}{4}$.
Then, $I^{N,p}(v)$ is attained at a function $Q$ satisfying the following.
\begin{itemize}
    \item[(i)] $\alpha = \frac{\sigma \|Q\|_2^{p-1}}{\omega}[\frac{4+N+p(4-N)}{N(p-1)\omega}]^{\frac{N}{4}-\frac{N(p+1)}{8}}[1+\frac{N}{2(p+1)}-\frac{N}{4}];$
 \item[(ii)] $Q$ is the ground state solution to Eq.~\eqref{eq64} with $\mu=0$,
  characterized by attaining the minimal $L^2$-norm of
 \begin{equation}\label{eq64.1}
\Delta^2 Q + \omega Q - \sigma |Q|^{p-1} Q = 0;
\end{equation}
  \item[(iii)] The  solution $Q$ of Eq.~\eqref{eq64.1}, that is, of \eqref{eq64} with $\mu=0$, satisfies the following Pohozaev identities:
\end{itemize}
\be
\int|\Delta Q|^{2}dx+\omega\int|Q|^{2}dx-\sigma\int|Q|^{p+1}dx=0,\nonumber
\ee
\be\label{eq9}
(2-\frac{N}{2})\int|\Delta Q|^2dx\:-\:\frac{\omega N}{2}\int|Q|^2dx+\frac{\sigma N}{p+1}\int|Q|^{p+1}dx\:=\:0.
\ee
\end{prop}

\begin{rem} In fact, from  Proposition \ref{prop2.4}, one can obtain the following sharp G-N inequality 
  \[
\|f\|_{p + 1}^{p + 1} \leqslant \frac{1}{\alpha} \|f\|_2^{p + 1 - \frac{p - 1}{4}N} \|\Delta f\|_2^{\frac{(p - 1)N}{4}},
\]
which was obtained by Fibich et al \cite{Fibich2002} and Zhu et al \cite{ZhuZhangYang2011} respectively. The existence of $Q$ was first obtained by
Zhu et al \cite{ZhuZhangYang2011}, and the uniqueness of minimal $L^{2}$-norm for Eq.~\eqref{eq64.1} was also proved in this paper. Recently, Ma \cite{Ma2023} proved the existence of positive solution for Eq.~\eqref{eq64.1} by the Schwarz symmetric method.
 Due to the presence of the biharmonic term, the corresponding fourth-order elliptic equation does not satisfy the classical maximum principle. As a result, traditional analytical tools for proving the uniqueness and positivity of solutions are no longer applicable. This necessitates the use of variational methods and profile decomposition techniques to establish the existence and structure of ground state solutions (see Zhu et al \cite{ZhuZhangYang2010}, Section 3 and Remark 3.3).
\end{rem}

\renewcommand{\theequation}{\thesection.\arabic{equation}}
\setcounter{equation}{0}
\section{New Variational Characterization}

For $u\in H^{2}$ and $\omega>0$, we define the following key functionals
\begin{align}
&T_{\omega,\mu}(u):=\int|\Delta u|^{2}dx+\mu\int|\nabla u|^{2}dx+\omega\int|u|^{2}dx,\\
 & F(u):=\frac{1}{2}T_{\omega,\mu}(u)-\frac{\sigma}{p+1}\int|u|^{p+1}dx,\label{eqI}\\
 & G(u):=T_{\omega,\mu}(u)-\sigma\int|u|^{p+1}dx,\label{eqS}\\
 & H(u):=8\int|\Delta u|^{2}dx+4\mu\int|\nabla u|^{2}dx-\frac{2N\sigma(p-1)}{p+1}\int|u|^{p+1}dx\label{eqQ}.
\end{align}
Consider the constrained variational problem
\begin{equation}\label{eqN}
  D_{\mathcal{M}^\mu_0}:=\inf\limits_{\mathcal{M}^\mu_0}F(u),
\end{equation}
where the  Nehari manifold $\mathcal{M}^\mu_0$ defined by 
$$\mathcal{M}^\mu_0:=\left\{u\in H^2\setminus\{0\},\quad G(u)=0\right\}$$
and we denote $\mathcal{M}^0_0:=\mathcal{M}^\mu_0|_{\mu=0}$. 
Next, we study the minimization problem \eqref{eqN} on the $\mathcal{M}^\mu_0$, and prove that $D_{\mathcal{M}^\mu_0}=\inf\limits_{\mathcal{M}^\mu_0}F(u)= \min\limits_{\mathcal{M}^\mu_0}F(u)$.

Here, we define an auxiliary functional 
$$J(u)=(\frac{1}{2}-\frac{1}{p+1})T_{\omega,\mu}(u).$$
Since $p>1$, it is easy to know that $J(u)$ is a non-negative functional, and from the expression of $J(u)$, one can deduce that
\begin{align}
F(u)&=\frac{1}{2}T_{\omega,\mu}(u)-\frac{\sigma}{p+1}\|u\|^{p+1}_{p+1}\nonumber\\
&=\frac{1}{2} G(u)+(\frac{1}{2}-\frac{1}{p+1})\sigma\|u\|^{p+1}_{p+1}\label{eq34}\\
&=\frac{1}{p+1}G(u)+J(u)\label{eq35}.
\end{align}
Then, $D_{\mathcal{M}^\mu_0}$ is rewritten as
\be\label{eq36}
D_{\mathcal{M}^\mu_0}=\inf_{\mathcal{M}^\mu_0}J(u).
\ee

\begin{lem}\label{lem3.1}
Let $N \geqslant  2$ and $\frac{N+8}{N}\leqslant  p<\frac{N+4}{ (N-4)^{+} }$. Then, $ D_{\mathcal{M}^\mu_0}\:>\:0$.
\end{lem}
\begin{proof} 
For all  $u \in \mathcal{M}^\mu_0$,  one can see 
\begin{equation*}\label{eq3}
T_{\omega,\mu}(u)=\sigma \int|u|^{p+1}dx, 
\end{equation*}
and 
\begin{equation*}
\begin{split}\label{eq4}
  F\left(u\right)=(\frac{1}{2}-\frac{1}{p+1})\sigma\int|u|^{p+1}dx=(\frac{1}{2}-\frac{1}{p+1})T_{\omega,\mu}(u).
\end{split}
\end{equation*}
It follows from the Sobolev embedding theorem that 
 $
  \sigma\int|u|^{p+1}dx \leqslant C\left(T_{\omega,\mu}(u)\right)^{\frac{p+1}2}$, 
where $C$ is positive constant.  Hence, one obtains
\begin{align*}
  T_{\omega,\mu}(u)\leqslant   CT^{\frac{p+1}2}_{\omega,\mu}(u),
\end{align*}
which implies 
$T_{\omega,\mu}(u) \geqslant  C^{\frac{2}{1-p}}>0$. And 
when $\frac{N+8}{N}\leqslant  p<\frac{N+4}{ (N-4)^{+} }$, one deduces that 
\begin{equation*}
F(u)=(\frac{1}{2}-\frac{1}{p+1})T_{\omega,\mu}(u) \: \geqslant  \: (\frac{1}{2}-\frac{1}{p+1})C^{\frac{2}{1-p}}:=C_3>0,
\end{equation*}
namely,
$ 
D_{\mathcal{M}^\mu_0}\geqslant  C_3 >0$. This completes the proof. 

\end{proof}

\begin{lem}\label{lem14}
Consider $N \geqslant 2$ and the exponent range $\frac{N+8}{N} \leqslant p < \frac{N+4}{(N-4)^+}$.
Let $\{u_n\} \subset H^2$ be a minimizing sequence for the variational problem \eqref{eqN} satisfying
\[
G(u_n) = 0 \quad \text{and} \quad F(u_n) \to D_{\mathcal{M}^\mu_0}.
\]
Then there exists a nonzero function $u \in H^2$ such that
\[
F(u) = D_{\mathcal{M}^\mu_0}.
\]

\end{lem}
\begin{proof}
Since $\{u_n\}$ is the minimizing sequence of the variational problem (\ref{eqN}), it follows from (\ref{eq36}) that $G(u_n)= 0$, and $J(u_n) \to D_{\mathcal{M}^\mu_0} $ as $n\to\infty$. From $\frac{N+8}{N}\leqslant  p<\frac{N+4}{ (N-4)^{+} }$, one sees that $\{u_n\}$ is bounded in $H^2$, and 
applying Lemma \ref{lemdi51}, $\{u_n\}$ can be decomposed as $u_n=\sum_{j = 1}^{l}U^j+u_n^l$ for every $l \geqslant 1$,  where $u_n^l = u_n^l(x)$ satisfies
\[
\lim_{l\to+\infty}\limsup_{n\to+\infty}\|u_n^l\|_q = 0,\quad\text{for every }q\in\left(2,\frac{2N}{(N-4)^+}\right),
\]
and (\ref{eq:1.8})-(\ref{eq:1.08}) holds. 
According to Remark 2.3 in \cite{ZhuZhangYang2011}, the following estimates hold as well.
\begin{align*}
\|u_n\|_{p+ 1}^{p + 1}&=\sum_{j = 1}^{l}\|U^j\|_{p+ 1}^{p + 1}+\|u_n^l\|_{p+ 1}^{p + 1}+o(1),\\
\|\nabla u_n\|_2^2&=\sum_{j = 1}^{l}\|\nabla U^j\|_2^2+\|\nabla u_n^l\|_2^2+o(1),
\end{align*}
where $o(1)=o_n(1)\to0$ as $n\to+\infty$. Injecting above estimates into $G(u)$ and $J(u)$, one deduces that 
\begin{align}
G(u_n)&=\sum_{j = 1}^{l}G(U^j)+G(u_n^l)+o(1),\label{eqdi2}\\
J(u_n)&=\sum_{j = 1}^{l}J(U^j)+J(u_n^l)+o(1)\label{eqdi1}.
\end{align}
It can be shown that $G(U^j) \geqslant 0$ for each $j \in \{1, \dots, l\}$. Suppose that, by contradiction, there exists $1 \leqslant j^* \leqslant l$ such that $G(U^{j^*}) < 0$. Then one can choose $\lambda \in (0,1)$ satisfying $G(\lambda U^{j^*}) = 0$, which implies that $\lambda U^{j^*} \in \mathcal{M}^\mu_0$, and hence $J(\lambda U^{j^*}) \geqslant D_{\mathcal{M}^\mu_0}$. 
Furthermore, it follows that
\begin{align*}
\lambda^2J(U^{j^*})=&\lambda^2\left[\left(\frac{1}{2}-\frac{1}{p + 1}\right)(\|\Delta U^{j^*}\|_2^2+\mu \|\nabla U^{j^*}\|_2^2+\omega\|U^{j^*}\|_2^2)\right]\\
=&J(\lambda U^{j^*})
\geqslant D_{\mathcal{M}^\mu_0}
\geqslant J(U^{j^*}),
\end{align*}
where the last inequality follows from \eqref{eqdi1} and the non-negativity of $J(u)$. 
This contradicts the fact that $0 < \lambda < 1$. 
Therefore, $G(U^j) \geqslant 0$ for all $1 \leqslant j \leqslant l$. 
In fact, $G(U^j) = 0$ follows directly from \eqref{eqdi2} together with the condition $G(u_n) = 0$.

Finally, the variational problem is attained. 
Indeed, there must exist $U^{\tilde{j}} \neq 0$ such that $G(U^{\tilde{j}}) = 0$ and $J(U^{\tilde{j}}) = D_{\mathcal{M}^\mu_0}$. 
Otherwise, if $U^{\tilde{j}} = 0$ for all $1 \leqslant j \leqslant l$, then this would contradict the fact that $J(U^{\tilde{j}}) = D_{\mathcal{M}^\mu_0} > 0$. 
Hence, there exists $U^{\tilde{j}} \neq 0$ with $G(U^{\tilde{j}}) = 0$, which implies $U^{\tilde{j}} \in \mathcal{M}^\mu_0$. Thus,
\begin{equation*}
J(U^{\tilde{j}})\leqslant D_{\mathcal{M}^\mu_0}\leqslant J(U^{\tilde{j}}),
\end{equation*}
where the first inequality can be deduced from \eqref{eqdi1} and the second follows from \eqref{eqN}. 
By setting $U^{\tilde{j}}=u$, one can obtain the equality $J(u)=F(u)=D_{\mathcal{M}^\mu_0}$, which establishes the desired result.
\end{proof}
Now we define cross-manifold
\be\label{eqM}
\mathcal{M}_-^0:=\left\{u\in H^2\setminus\{0\}, \,\ G(u)<0, \:H(u)=0\right\},
\ee
and focus on the minimization problem with cross constraints
\be\label{eq54}
D_{\mathcal{M}_-^0}:=\inf_{\mathcal{M}_-^0}F(u).
\ee

\begin{lem}\label{lem3.2} Let $\mu=0$, $N\geqslant 2$ and $\frac{N+8}{N}\leqslant  p<\frac{N+4}{ (N-4)^{+} }$. 
There exists $u \in H^2 \setminus \{ 0\}$ such that $G(u)=0$ and $H(u)=0$.
\end{lem}
\begin{proof}
Observe that if $u(x)$ is a bound state of the elliptic equation
\be\label{eq8} 
\Delta^2 u+\omega u -\sigma|u|^{p-1}u=0.
\ee
Then, from the Pohozaev identities in Proposition \ref{prop2.4}, a direct conclusion can be drawn.

\end{proof}

\begin{lem}\label{lem3.3}
Let $\mu=0$ and $N \geqslant 2$. Then the set $\mathcal{M}_-^0$ is nonempty for all $p$ satisfying $\frac{N+8}{N}\leqslant  p<\frac{N+4}{ (N-4)^{+} }$.
\end{lem}
\begin{proof}
From Lemma \ref{lem3.2}, there exists $u\in H^2\setminus\{0\}$ such that both $G(u)=0$ and $H(u)=0$. 
Let $u_\beta=\beta^{\frac{4}{p-1}} u(\beta x)$ for $\beta >0$. Then,  $u_\beta \in H^2\setminus\{0\}$, by \eqref{eqI} and \eqref{eqS}, one gets
\begin{equation*}\label{eq10}
  \begin{split}
  G(u_{\beta})&=\beta^{\frac{8}{p-1}+4-N}(\int|\Delta u|^2dx-\sigma\int|u|^{p+1}dx)+\beta^{\frac{8}{p-1}-N}\omega\int|u|^2dx\\
 &=\beta^{\frac{8}{p-1}-N}(1-\beta^4)\:\omega\:\int|u|^2dx.
  \end{split}
\end{equation*} 
Thus, $G(u)=0$ implies that there exists $\beta^*>1$ such that $G(u_{\beta^*})<0$. On the other hand, one deduces that 
\begin{equation*}\label{eq11}
\begin{split}
H(u_\beta)&=8\beta^{\frac{8}{p-1}+4-N}\int|\Delta u|^{2}dx-\frac{2N\sigma(p-1)}{p+1} \beta^{\frac{4(p+1)}{p-1}-N}\int|u|^{p+1}dx\\
&=\beta^{\frac{8}{p-1}+4-N}H(u).
\end{split}
\end{equation*} 
This implies that $H(u_{\beta^*})=(\beta^*)^{\frac{8}{p-1}+4-N}H(u)=0$, so
$u_{\beta^*}\in\mathcal{M}_-^0$.

\end{proof}
\begin{lem}\label{lem3.4} 
  Let $\mu=0$ and $N \geqslant  2$. Then $D_{\mathcal{M}_-^0} > 0$ provided $\frac{N+8}{N}\leqslant  p<\frac{N+4}{ (N-4)^{+} }$.
\end{lem}
\begin{proof}
Let $u\in\mathcal{M}_-^0$, so that $u\neq0$, $G(u)<0$ and $H(u)=0$. Then
\be\label{eq12}
\begin{split}
  F\left(u\right)&=\frac{N(p-1)-8}{8(p+1)}\sigma\int|u|^{p+1}dx\:+\:\frac \omega2 \int|u|^2dx.
\end{split}
\ee
Since $p\geqslant\frac{N+8}{N}$, both terms on the right-hand side of \eqref{eq12} are nonnegative, whence $D_{\mathcal{M}_-^0}\geqslant 0$. In the following, the conclusion $D_{\mathcal{M}_-^0}> 0$ is discussed in two cases.

Case $1)$: The critical case $p= \frac{N+8}{N}$;

Case $2)$: The supercritical case $\frac{N+8}{N}< p<\frac{N+4}{(N-4)^{+}}$.
 
To begin with Case 1), noting that $D_{\mathcal{M}_-^0} \geqslant 0$, the proof proceeds via contradiction. Assume that $D_{\mathcal{M}_-^0} = 0$. Then, from \eqref{eq54}, one can find a sequence $\{u_n\}_{n=1}^{\infty}\subset \mathcal{M}_-^0$ for $n \in \mathbb{Z}^+$ satisfying
\[
H(u_n) = 0, \quad G(u_n) < 0, \quad \text{and} \quad F(u_n) \to 0 \quad \quad \text{as } n \to \infty.
\]
Due to \eqref{eq12} and $p=\frac{N+8}{N}$, one gets
\be\label{eq52}
\int|u_n|^2dx\to0\quad \quad \text{as}\,\ n\to\infty.
\ee
On the one hand, it follows from G-N inequality and $p=\frac{N+8}{N}$ that 
\be\label{eq50}
\int|u_n|^{p+1}dx \lesssim (\int|\Delta u_n|^2dx)\cdot\left(\int|u_n|^2dx\right)^{\frac4N},
\ee
according to $G(u_n)<0$ and \eqref{eq50}, there exists a constant $C>0$, independent of $n$,
such that
\begin{align}
  \int|\Delta u_{n}|^{2}dx+\omega\int|u_{n}|^{2}dx&  \leqslant C\Big(\int|\Delta u_n|^2dx\Big)\Big(\int|u_n|^2dx\Big)^{\frac4N}\label{eq51}.
\end{align}
Since $u_n\in\mathcal{M}_-^0$, we have $u_n\neq0$ and therefore $\int|\Delta u_n|^2dx>0$;
dividing \eqref{eq51} by this quantity gives
\be\label{eq53}
1\leqslant\frac{\int|\Delta u_n|^2dx+\omega\int|u_n|^2dx}{\int|\Delta u_n|^2dx}
 \leqslant C\Big(\int|u_n|^2dx\Big)^{\frac4N}.
\ee
By \eqref{eq52}, the right-hand side of \eqref{eq53} tends to $0$ as $n\to\infty$, which is a
contradiction. Therefore $D_{\mathcal{M}_-^0}>0$ for $p=\frac{N+8}{N}$.

Next, for Case $2)$,
from the Sobolev embedding inequality and Eq.~\eqref{eq76}, one  deduce that
\begin{equation*}
\begin{split}
  \int|u|^{p+1}dx& \lesssim \left(\int|\nabla u|^2dx+\int|\Delta u|^2dx+\omega\int|u|^2dx\right)^{\frac{p+1}2}\\
&\leqslant  C\,T_{\omega,0}(u)^{\frac{p+1}2},
\end{split}
\end{equation*}
where the second inequality uses $\int|\nabla u|^2dx\leqslant\|u\|_{2}\|\Delta u\|_{2}$ and
$\omega>0$. According to $G(u)<0$, that is $T_{\omega,0}(u)<\sigma\int|u|^{p+1}dx$, it follows that 
\be 
  T_{\omega,0}(u) \leqslant   \sigma C \,T_{\omega,0}(u)^{\frac{p+1}2}.
\ee
Since $\frac{N+8}{N}< p<\frac{N+4}{ (N-4)^{+} }$ and $N \geqslant  2$, one deduces that
\be \label{eq14}
T_{\omega,0}(u)\geqslant  (\sigma C)^{\frac{2}{1-p}}>0.
\ee
Thus, by \eqref{eq12}, the following inequality holds
\begin{equation*}
\begin{split}
  F\left(u\right)&> \frac{N(p-1)-8}{8(p+1)}\,T_{\omega,0}(u)\:+\:\frac \omega2 \int|u|^2dx\\
  &\geqslant   \frac{N(p-1)-8}{8(p+1)}(\sigma C)^{\frac{2}{1-p}}>0.
\end{split}
\end{equation*} 
Here, the first inequality follows from the condition $G(u) < 0$, while the second one is justified by \eqref{eq14} and $N(p-1)>8$.
 Hence, $D_{\mathcal{M}_-^0}=\inf\limits_{\mathcal{M}_-^0}F(u)>0$ for the case $\frac{N+8}{N}< p<\frac{N+4}{(N-4)^{+}}$.

In conclusion, the proof of Lemma \ref{lem3.4} is complete.

\end{proof}
Now, we state the main theorem of this section. More precisely, we obtain a comparison between  $D_{\mathcal{M}^0_0}$ and $D_{\mathcal{M}_-^0}$, which provides a new variational characterization to investigate the sharp threshold of blow-up and global existence for the time-dependent biharmonic nonlinear Schr\"odinger equation. This result is also crucial for proving the strong instability of standing waves without any frequency constraint. 
 
\begin{thm}\label{lem12}
Let $\mu=0$, $N\geqslant 2$ and $\frac{N+8}{N}\leqslant  p<\frac{N+4}{ (N-4)^{+} }$. 
Then, one has  $D_{\mathcal{M}^0_0}\leqslant  D_{\mathcal{M}_-^0}$.
\end{thm}

\begin{proof}
We distinguish two cases according to the exponent $p$.

Case $1)$: The critical case $p= \frac{N+8}{N}$;

Case $2)$: The supercritical case $\frac{N+8}{N}< p<\frac{N+4}{(N-4)^{+}}$.
 
For Case $1)$, 
let
$v\in\mathcal{M}_-^0$. Then $G(v)<0$ and $H(v)=0$. Since $N(p-1)=8$, the condition $H(v)=0$
gives $\frac{\sigma}{p+1}\int |v|^{p+1}dx=\frac12\int |\Delta v|^2dx $.
For the scaling $v_\eta(x)=\eta^{\frac N4}v(\eta^{\frac12}x)$ and $\eta>0$, one deduces that
\begin{equation*}
  F(v_\eta)=\frac{\omega}{2}\int |v|^2dx
  +\eta^2\left(\frac12\int |\Delta v|^2dx-\frac{\sigma}{p+1}\int |v|^{p+1}dx\right)
  =\frac{\omega}{2}\int |v|^2dx=F(v).
\end{equation*}
On the other hand, by \eqref{eqS} and $\frac{N(p-1)}{4}=2$,
\begin{equation*}
  G(v_\eta)=\omega\int |v|^2dx
  +\eta^2\left(\int |\Delta v|^2dx-\sigma\int |v|^{p+1}dx\right).
\end{equation*}
Since $G(v)=G(v_1)<0$, while
$\displaystyle  \lim_{\eta\to0^+}G(v_\eta)=\omega\int |v|^2dx>0$, the continuity of
$\eta\mapsto G(v_\eta)$ implies that there exists $\eta^*\in(0,1)$ such that
$G(v_{\eta^*})=0$. Thus $v_{\eta^*}\in\mathcal{M}^0_0$, and consequently
$
  D_{\mathcal{M}^0_0}\leqslant F(v_{\eta^*})=F(v)$. 
Taking the infimum over all $v\in\mathcal{M}_-^0$ yields
\begin{equation*}
  D_{\mathcal{M}^0_0} \leqslant \inf\limits_{v\in \mathcal{M}_-^0 }F(v)= D_{\mathcal{M}_-^0}.
\end{equation*}

For Case $2)$:  $\frac{N+8}{N}< p<\frac{N+4}{ (N-4)^{+} }$, 
let  $u \in H^2\setminus\{0\}$ with $H(u)\leqslant 0$. By using the scaling transformation
$u_\eta (x)=\eta^{\frac{N}{4}}u(\eta^{\frac{1}{2}}x)$ for $\eta>0$, one deduces that
\be\label{eq66}
  F\left(u_\eta\right)=\frac{\omega}{2}\int|u|^{2}dx+\frac{\eta^{2}}{2}\int|\Delta u|^{2}dx
  -\frac{\eta^{\frac{N(p-1)}{4}}\sigma}{p+1}\int|u|^{p+1}dx ,
\ee
\begin{align}\label{eq67}
  H\left(u_\eta\right)=&\eta^{2}\Big(8\int|\Delta u|^{2}dx
  -\eta^{\frac{N(p-1)-8}{4}}\frac{2N\sigma(p-1)}{p+1}\int|u|^{p+1}dx\Big)\\
  =:&\eta^{2}h\left(\eta\right)\nonumber.
\end{align}
Since $p>\frac{N+8}{N}$, one has $\frac{N(p-1)-8}{4}>0$, so that $h\left(\eta\right)$ is
strictly decreasing on $(0,+\infty)$ with
$\displaystyle \lim_{\eta\to0^+}h(\eta)=8\int|\Delta u|^{2}dx>0$. Hence, there exists a unique $\eta_0 \in (0,1]$ such
that $H\left(u_{\eta_0}\right)=0$ and
\be\label{eq83}
  H(u)=0 \quad\text{if and only if}\quad \eta_0=1
\ee
by \eqref{eq67} and  $H(u)\leqslant 0$. Moreover, differentiating \eqref{eq66} with respect to $\eta$, one obtains
\begin{equation*}
  \frac{d}{d\eta}F\left(u_\eta\right)=\frac{H\left(u_\eta\right)}{8\eta} .
\end{equation*}
From this and the fact that $H\left(u_{\eta_0}\right)=0$, it can be easily checked that
\[
\frac{d}{d\eta} F(u_{\eta}) =
\begin{cases}
> 0, & \text{if }\eta \in (0, \eta_0), \\
< 0, & \text{if } \eta \in (\eta_0, \infty).
\end{cases}
\]
It follows from the shape of the graph of the functional $F(u_\eta)$ that
\be\label{eq69}
F(u_{\eta})<F(u_{\eta_0}), \quad \quad\text{for any}\,\ \eta>0, \,\ \eta \neq \eta_0.
\ee

Now, the
relationship between $D_{\mathcal{M}^0_0}$ and $D_{\mathcal{M}_-^0}$ can be proved. Let
$v\in \mathcal{M}_-^0$, that is, $G(v)<0$ and $H(v)=0$, so that $\eta_0=1$ by \eqref{eq83}.
Scaling $v_\eta (x)=\eta^{\frac{N}{4}}v(\eta^{\frac{1}{2}}x)$ yields
\be\label{eq68}
  G\left(v_\eta\right)=\omega\int| v|^{2}dx+\eta^{2}\int|\Delta v|^{2}dx
  -\eta^{\frac{N(p-1)}{4}}\sigma\int|v|^{p+1}dx .
\ee
Noted that $G(v_\eta)>0$ for $\eta$
small enough. Thus, there exists a parameter $\eta^* \in (0,1)$ such that
$G(v_{\eta^*})=0$, which implies that $v_{\eta^*} \in \mathcal{M}^0_0$ and
\begin{equation*}
  D_{\mathcal{M}^0_0} \leqslant F(v_{\eta^*})< F(v),
\end{equation*} 
where the second inequality follows from \eqref{eq69} and $\eta_0=1$. Taking the infimum of
$F(v)$ in \eqref{eq54} yields the desired conclusion.

Combining Cases 1) and 2) yields $D_{\mathcal{M}^0_0}\leqslant D_{\mathcal{M}_-^0}$,
which completes the proof.

  \end{proof}

\renewcommand{\theequation}{\thesection.\arabic{equation}}
\setcounter{equation}{0}
\section{Sharp Criteria of blow-up}

We analyze in this section the sharp conditions that govern blow-up and global existence of solutions to Eq.~\eqref{eq1} under the cases $\mu=0$ and $\mu>0$.
 It can reveal the contribution of the lower-order dispersive term to blow-up solutions in Eq.~\eqref{eq1}. To begin with, we investigate the sharp criteria for blow-up solutions of Eq.~\eqref{eq1} by utilizing the cross-constraint variational arguments.
\subsection{Sharp Criterion in the case $\mu$=0}
We introduce the invariant manifold as follows
\be\label{eqK}
\mathcal{W}^{0}_{b}:=\left\{\phi\in H^2,\quad F(\phi)<D_{\mathcal{M}^0_0},\quad G(\phi)<0,\quad H(\phi)<0\right\}.
\ee
Thus, the following propositions hold.
\begin{prop}\label{prop3.71}
Let $\mu=0$, $N \geqslant 2$ and $\frac{N+8}{N} \leqslant p < \frac{N+4}{(N-4)^+}$. Then the set $\mathcal{W}^{0}_{b}$ is nonempty and forms an invariant manifold for the Cauchy problem \eqref{eq1}. More precisely, if the initial data $\phi_0(x)$ belongs to $\mathcal{W}^{0}_{b}$, then the corresponding solution $\phi(t,x)$ remains in $\mathcal{W}^{0}_{b}$ for all $t \in [0,T)$.
\end{prop}
\begin{proof} 
First of all, it is shown that $\mathcal{W}^{0}_{b}$ is nonempty.
From Lemma \ref{lem3.2}, there exists $u\in  H^2\setminus\{0\}$ such that both $G(u)=0$ and
$H(u)=0$.  It follows from \eqref{eqI}, \eqref{eqS} and \eqref{eqQ} that 
\begin{align}
  F(\lambda u)\:=&\:\frac{1}{2}\lambda^{2}\left(\int|\Delta u|^{2}dx+\omega\int|u|^{2}dx\right) -\frac{\sigma}{p+1}\lambda^{p+1}\int|u|^{p+1}dx,\label{eq16.1}\\
  G(\lambda u)\:=&\:G(\lambda u)-\lambda^{2}G(u)\:=(1-\lambda^{p-1})\lambda^{2}\sigma \int|u|^{p+1}dx, \label{eq16.2}\\
  H(\lambda u)\:=&\:H(\lambda u)-\lambda^{2}H(u)\:=(1-\lambda^{p-1})\lambda^{2}\frac{2N\sigma(p-1)}{p+1} \int|u|^{p+1}dx.\label{eq16.3}
\end{align}
$D_{\mathcal{M}^0_0}>0$ and $\frac{N+8}{N}\leqslant  p<\frac{N+4}{ (N-4)^{+} }$ implies that, for $\lambda>1$ large enough,   
\[G(\lambda u)<0, \,\ H(\lambda u)<0\ \ {\rm and }\ \ F(\lambda u)<D_{\mathcal{M}^0_0}.\] It means  $\lambda u\in\mathcal{W}^{0}_{b}$. Namely $\mathcal{W}^{0}_{b}$ is not empty.

Then, assuming $\phi_0 \in \mathcal{W}^{0}_{b}$, Lemma \ref{prop2.1} guarantees the existence of a unique solution $\phi(t, \cdot) \in C([0,T); H^2)$ 
to the Cauchy problem \eqref{eq1} with initial data $\phi_0$. By the identities \eqref{eqm}, \eqref{eqe1}, and \eqref{eqI}, the functional $F$ is conserved with respect to time $t$:
\begin{equation}\label{eq17.2}
F(\phi(t,\cdot)) = F(\phi_0), \quad t \in [0,T).
\end{equation}
Hence, the inequality $F(\phi_0) < D_{\mathcal{M}^0_0}$ implies that this bound holds for all $t\in[0,T)$.

Next, it is proved by contradiction that $G(\phi(t,\cdot)) < 0$ for all $t \in [0,T)$.
 Suppose otherwise, then by continuity, there exists $t^* \in (0,T)$ such that $G(\phi(t^*,\cdot)) = 0$. From \eqref{eq17.2}, it follows that $\phi(t^*,\cdot) \not\equiv 0$. By \eqref{eqN} and Theorem \ref{lem12}, one gets $F(\phi(t^*,\cdot)) \geqslant D_{\mathcal{M}^0_0}$,
which contradicts the fact that $F(\phi(t,\cdot)) < D_{\mathcal{M}^0_0}$ for all $t \in [0,T)$. Hence, $G(\phi(t,\cdot)) < 0$ holds for every $t \in [0,T)$.

Finally, the claim $H(\phi(t,\cdot)) < 0$ is established by contradiction. If not, there exists $\widetilde{t} \in (0,T)$ such that $H(\phi(\widetilde{t},\cdot)) = 0$. Since it has been shown above that $G(\phi(\widetilde{t},\cdot)) < 0$, it follows that $\phi(\widetilde{t},\cdot) \in \mathcal{M}_-^0$. Then, \eqref{eqM} and Theorem \ref{lem12} deduce that
$F(\phi(\widetilde{t},\cdot)) \geqslant D_{\mathcal{M}_-^0} \geqslant D_{\mathcal{M}^0_0}$,
which contradicts the assumption $F(\phi(t,\cdot)) < D_{\mathcal{M}^0_0}$ for all $t \in [0,T)$. Therefore, $H(\phi(t,\cdot)) < 0$ holds for every $t \in [0,T)$. This completes the proof of Proposition \ref{prop3.71}.

\end{proof}
Using an argument analogous to that in Proposition \ref{prop3.71}, the following invariant manifolds can be established.
\begin{prop}\label{prop3.8}
  Let $\mu=0$, $N \geqslant  2$ and $\frac{N+8}{N}\leqslant  p <\frac{N+4}{ (N-4)^{+} }$. Define
\begin{align*}
  \mathcal{W}^{0}_{g_1}:&=\left\{\phi\in H^{2},\quad F(\phi)<D_{\mathcal{M}^0_0},\quad G(\phi)<0,\quad H(\phi)>0\right\},\\
  \mathcal{W}_{0}:&=\left\{\phi\in H^{2},\quad F(\phi)<D_{\mathcal{M}^0_0},\quad G(\phi)<0\right\},\\
  \mathcal{W}^{0}_{g_2}:&=\left\{\phi\in H^{2},\quad F(\phi)<D_{\mathcal{M}^0_0},\quad G(\phi)>0\right\}.
\end{align*}
Consequently, $\mathcal{W}^{0}_{g_1}$, $\mathcal{W}_{0}$, and $\mathcal{W}^{0}_{g_2}$ serve as invariant manifolds for the Cauchy problem \eqref{eq1}.
\end{prop}

\begin{rem}\label{rem3.9}
Among the manifolds introduced in Proposition \ref{prop3.8}, $\mathcal{W}_{0}$ and $\mathcal{W}^{0}_{g_2}$ are referred to  invariant manifolds for the Cauchy problem \eqref{eq1}. Naturally, the sets $\mathcal{W}^{0}_{b}$ and $\mathcal{W}^{0}_{g_1}$ are termed cross-invariant manifolds related to the same problem.
\end{rem}
\vspace{0.2cm}
Now, the sharp sufficient conditions for blow-up and global existence are discussed. To begin with, a sufficient criterion for the blow-up solution to the Cauchy problem \eqref{eq1} is presented.

\vspace{0.3cm}
{\bf Proof of Theorem \ref{th1}-{i}}: 
The proof proceeds in three steps.

Step 1: It is claimed that there exists a constant 
$$a:= D_{\mathcal{M}^0_0}-F(\phi_0) > 0,$$ 
such that $H(\phi) \leqslant -16a$ for all time $t$.
According to $\phi_0\in\mathcal{W}^{0}_{b}$ and Proposition \ref{prop3.71}, which implies that $\phi(t,.)\in\mathcal{W}^{0}_{b}$ is true for $t\in[0,T)$. Thus, one deduces that $F(\phi)<D_{\mathcal{M}^0_0},\,\ G(\phi)<0$ and $H(\phi)<0$ are true for all time $t$. 
Consider the scaling  $\phi_\lambda=\lambda^{\frac N2}\phi(\lambda x)$. 
\begin{align}\label{eq20}
&G(\phi_{\lambda})=\omega\int|\phi|^{2}dx+\lambda^{4}\int|\Delta \phi|^{2}dx- \lambda^{\frac{N(p-1)}{2}}\sigma\int|\phi|^{p+1}dx, \\
&H(\phi_{\lambda})=8\lambda^{4}\int|\Delta \phi|^{2}dx-\lambda^{\frac{N(p-1)}{2}}\frac{2N\sigma(p-1)}{p+1}\int|\phi|^{p+1}dx.
\end{align}
Since $G(\phi) < 0$, there exists a parameter $\lambda^* \in (0,1)$ such that
$G(\phi_{\lambda^*}) = 0$
and moreover, 
$G(\phi_{\lambda}) < 0 \,\ \text{for all } \lambda \in (\lambda^*, 1]$.
It follows from $H(\phi) < 0$ that the functional $H(\phi_{\lambda})$ exhibits one of the following three scenarios:

1-i) $H( \phi _{\lambda }) < 0$ for every $\lambda \in [ \lambda ^{* }, 1];$\par
 1-ii) $H( \phi _{\lambda }* ) = 0$;\par
 1-iii) There exist $\xi  \in ( \lambda ^{* }, 1)$ such that $H( \phi _{\xi }) = 0$. 

 For the first two cases, since $G(\phi_{\lambda^*}) = 0$ and $H(\phi_{\lambda^*}) \leqslant 0$, by \eqref{eqN}, \eqref{eqM} and Theorem \ref{lem12}, one obtains
$F(\phi_{\lambda^*}) \geqslant D_{\mathcal{M}^0_0}$.
Furthermore, it is straightforward to check that
\begin{align}
  &F(\phi)-F(\phi_{\lambda^{*}}) =\frac{1}{2}(1-\lambda^{*4})\int|\Delta \phi|^{2}dx -(1-\lambda^{*\frac{N(p-1)}{2}})\frac{\sigma}{p+1}\int|\phi|^{p+1}dx,\label{eq21}\\
  &H(\phi)-H(\phi_{\lambda^{*}}) =8(1-\lambda^{*4})\int|\Delta \phi|^{2}dx -(1-\lambda^{*\frac{N(p-1)}{2}})\frac{2\sigma N(p-1)}{p+1}\int|\phi|^{p+1}dx.\label{eq22}
\end{align}
By $0<\lambda^*<1$ and  $\frac{N+8}{N}\leqslant  p <\frac{N+4}{ (N-4)^{+} }$, it follows from \eqref{eq21} and \eqref{eq22}   that
\begin{align}\label{eq23}
F(\phi)-F(\phi_{\lambda^*})&\geqslant  \frac{1}{2}(1-\lambda^{*4})\int|\Delta \phi|^{2}dx-(1-\lambda^{*\frac{N(p-1)}{2}})\frac{\sigma N(p-1)}{8(p+1)}\int|\phi|^{p+1}dx \nonumber\\
&=\frac{1}{16}H(\phi)-\frac{1}{16}H(\phi_{\lambda^*})\geqslant \frac{1}{16}H(\phi).
\end{align}

For Case 1-iii), one sees  that  $H(\phi_\xi)=0$ and $G(\phi_\xi)<0$, which implies  $\phi_\xi\in \mathcal{M}_-^0$. From \eqref{eqM} and Theorem \ref{lem12}, it follows that $F(\phi_\xi)\geqslant  D_{\mathcal{M}_-^0}\geqslant  D_{\mathcal{M}^0_0}$. In addition,
\be\label{eq24}
F\left(\phi\right)-F\left(\phi_{\xi}\right)\geqslant \frac{1}{16}H(\phi)-\frac{1}{16}H(\phi_{\xi})\geqslant \frac{1}{16}H(\phi).
\ee
Since $F(\phi_{\lambda^*})\geqslant  D_{\mathcal{M}^0_0}$ and $F(\phi_\xi)\geqslant  D_{\mathcal{M}^0_0}$, in view of \eqref{eq23} and \eqref{eq24}, one gets
\be \label{eq25}
\begin{cases}H(\phi)\leqslant 16(F(\phi)-F(\phi_{\xi}))\leqslant 16[F(\phi)-D_{\mathcal{M}^0_0}],\\ H(\phi)\leqslant 16(F(\phi)-F(\phi_{\lambda^{*}}))\leqslant 16[F(\phi)-D_{\mathcal{M}^0_0}].
\end{cases}
\ee
From \eqref{eqm}, \eqref{eqe1} and \eqref{eqI}, it has
$$F\left(\phi\right)=F\left(\phi_0\right).$$
Thus by $\phi_0\in\mathcal{W}^{0}_{b}$ and \eqref{eq25}, one can show that easily
\be\label{eq70}
H(\phi) \leqslant  16[F(\phi_0)-D_{\mathcal{M}^0_0}]=-16a< 0.
\ee
where $a=D_{\mathcal{M}^0_0}-F(\phi_0) > 0$ is defined as above.

Step 2: There exists a constant $\tau >0$ such that
$$\frac{d}{dt} \mathcal{M}_{\varphi_R}[\phi(t)] \leqslant -\tau \|\Delta \phi(t)\|^{2}_{2}, \quad \text{for all} \,\ t\in [0, T),$$
and  $\mathcal{M}_{\varphi_R}[\phi(t)] < 0$ holds  for all $t \geqslant t_1$, where    $t_1 \geqslant 0$ is a fixed time. 
For this purpose, the analysis now is divided into the following two cases.

Case $1)$: Suppose that 
\begin{equation*}
\|\Delta \phi(t)\|^2_2 \leqslant \frac{4N(p-1) E(\phi_0)}{N(p-1) - 8}, \quad \text{for all} \,\ t\in [0, T).
\end{equation*}
Since \( \|\nabla\phi(t)\|^2_2 \leqslant \|\phi(t)\|_2\|\Delta\phi(t)\|_2 =\|\phi(0)\|_2\|\Delta\phi(t)\|_2 \), by taking \( R > 0 \) large enough, one derives that 
\begin{align*}
& O\left( R^{-4} + \|\nabla\phi(t)\|^2_2 R^{-2} + \|\nabla\phi(t)\|^{\frac{p-1}{2}}_2 R^{-\frac{(p-1)(N-1)}{2}} \right) \\
\lesssim  & O\left( R^{-4} + \|\Delta \phi(t)\|_2 R^{-2} + \|\Delta\phi(t)\|^{\frac{p-1}{4}}_2 R^{-\frac{(p-1)(N-1)}{2}} \right) \\
< &a.
\end{align*}
Combining Lemma \ref{lem11} with \eqref{eq70}, one obtains that
\begin{equation}\label{eq71}
\frac{d}{dt}\mathcal{M}_{\varphi_R}[\phi(t)] \leqslant -15a \leqslant -\tau\|\Delta \phi(t)\|^2_2
\end{equation}
for sufficiently small positive parameter \( \tau \).

Case $2)$: Assume that for all \( t \in [0, T) \) satisfies
\begin{equation*}
\|\Delta\phi(t)\|^2_2 > \frac{4N(p-1) E(\phi_0)}{N(p-1) - 8}.
\end{equation*}
Due to Lemma \ref{lem11}, then
\begin{align*}
\frac{d}{dt} \mathcal{M}_{\varphi_R}[\phi(t)] & \leqslant 2N(p-1) E(\phi_0) - (N(p-1) - 8)\|\Delta \phi(t)\|_{2}^2\\
&\quad +O\left(R^{-4} + R^{-2}\|\nabla \phi(t)\|_{2}^2 + R^{-\frac{(p-1)(N-1)}{2}}\|\nabla \phi(t)\|_{2}^{\frac{p-1}{2}}\right)\\
&\lesssim   -\frac{1}{2}(N(p-1) - 8)\|\Delta \phi(t)\|_{2}^2\\
&\quad + O\left( R^{-4} + R^{-2}\|\Delta \phi(t)\|_2 + R^{-\frac{(p-1)(N-1)}{2}}\|\Delta\phi(t)\|^{\frac{p-1}{4}}_2 \right).
\end{align*}
Note that $p\leqslant 9$. Then, take \( R \) large enough, the following estimate holds for all time $t$
\begin{equation}\label{eq72}
\frac{d}{dt}\mathcal{M}_{\varphi_R}[\phi(t)] \leqslant -\frac{1}{4}(N(p-1) - 8)\|\Delta \phi(t)\|_{2}^2.
\end{equation}
By combining the estimates \eqref{eq71} and \eqref{eq72}, one obtains the existence of a constant \( \tau > 0 \) such that
\be\label{eq73.1}
\frac{d}{dt} \mathcal{M}_{\varphi_R}[\phi(t)] \leqslant -\tau \|\Delta \phi(t)\|^{2}_{2}<0
\ee
holds for all $t\in [0, T)$. According to the strict monotonically decreasing function $\mathcal{M}_{\varphi_R}[\phi]$, there exists a $t_1 \geqslant 0$ such that 
\be\label{eq74.1}
\mathcal{M}_{\varphi_R}[\phi(t)] < 0, \quad \text{for} \,\ t \geqslant t_1.
\ee

Step 3: The solution \(\phi(t,x)\) of the Cauchy problem \eqref{eq1} blows up in  a  finite time.
Integrating \eqref{eq73.1} on $[t_1, t]$ and using \eqref{eq74.1}, one can deduce that 
\begin{align*}
\mathcal{M}_{\varphi_R}[\phi(t)] =\mathcal{M}_{\varphi_R}[\phi(t_1)] + \int_{t_1}^{t}\frac{d}{ds}\mathcal{M}_{\varphi_R}[\phi(s)]ds \leqslant  -\tau \int_{t_1}^{t}\|\Delta \phi(s)\|^{2}_{2}ds.
\end{align*}
On the other hand, by applying the Cauchy-Schwarz inequality, it follows from the definition of \(\mathcal{M}_{\varphi_R}[\phi(t)]\) in \eqref{eq75} that
\[
|\mathcal{M}_{\varphi_R}[\phi(t)]| \leqslant 2\|\nabla\varphi_R\|_\infty\|\phi(t)\|_2\|\nabla\phi(t)\|_2 \leqslant C(\phi_0)R\|\Delta\phi(t)\|^{\frac{1}{2}}_2.
\]
Namely,
\[
|\mathcal{M}_{\varphi_R}[\phi(t)]|^4  \leqslant C^4(\phi_0)R^4\|\Delta\phi(t)\|^{2}_2,
\]
and 
$$\int_{t_1}^t |\mathcal{M}_{\varphi_R}[\phi(s)]|^4ds \leqslant C^4(\phi_0)R^4 \int_{t_1}^t\|\Delta\phi(s)\|^{2}_2ds \leqslant -\frac{C^4(\phi_0)R^4}{\tau}\mathcal{M}_{\varphi_R}[\phi(t)].$$
Thus for $A:=\frac{\tau}{C^4(\phi_0)R^4}$, 
\begin{equation}\label{3.10}
\mathcal{M}_{\varphi_R}[\phi(t)] \leqslant -A\int_{t_1}^t |\mathcal{M}_{\varphi_R}[\phi(s)]|^4ds.
\end{equation}
Let $\gamma_1 (t):=\int_{t_1}^t |\mathcal{M}_{\varphi_R}[\phi(s)]|^4ds$ for $t\in[t_1,T)$. This
function is nonnegative and strictly increasing, since
$\mathcal{M}_{\varphi_R}[\phi(s)]<0$ for every $s\geqslant t_1$. Moreover, \eqref{3.10} gives
$
\gamma_1 '(t)=\big|\mathcal{M}_{\varphi_R}[\phi(t)]\big|^4 \geqslant A^4\gamma_1 (t)^4 .
$
Hence, fixing a time $t_2\in(t_1,T)$, so that $\gamma_1(t_2)>0$, and integrating this
differential inequality on $[t_2,t]$, one gets
\[
\gamma_1 (t) \geqslant \frac{\gamma_1 (t_2)}{(1 - 3A^4(t - t_2)\gamma_1 ^3(t_2))^{\frac{1}{3}}},
\]
and consequently, by \eqref{3.10},
\[
\mathcal{M}_{\varphi_R}[\phi(t)] \to -\infty \quad \quad\text{as} \quad \quad t \to t_1^* := t_2 + \frac{1}{3 A^4 \gamma_1^3(t_2)} < \infty.
\]
Since $|\mathcal{M}_{\varphi_R}[\phi(t)]| \leqslant C(\phi_0)R\|\Delta\phi(t)\|^{1/2}_2$, this
forces $\|\Delta\phi(t)\|_2\to+\infty$, and therefore $T\leqslant t_1^*<+\infty$.
This implies that the solution \(\phi\) of the Cauchy problem \eqref{eq1} blows up in a finite time under the assumptions \(\mu = 0\) and 
$\frac{N+8}{N} < p < \frac{N+4}{(N-4)^+}$.

\hfill \(\square\)

\vspace{0.2cm}
In what follows, a sufficient criterion for the global existence of solution to the Cauchy problem \eqref{eq1} will be provided.

\vspace{0.3cm}
{\bf Proof of Theorem \ref{th1}-{ii}}: 
This theorem is proved by two cases.

Case\:1): $\phi_0\in\mathcal{W}^{0}_{g_1}$;

Case\:2): $\phi_0\in\mathcal{W}^{0}_{g_2}$.

Firstly, consider the Case 1): \(\phi_0 \in \mathcal{W}^{0}_{g_1}\). By Proposition \ref{prop3.8}, the solution \(\phi(t,x)\) to the Cauchy problem \eqref{eq1} satisfies \(\phi(t,\cdot) \in \mathcal{W}^{0}_{g_1}\) for all \(t \in [0,T)\).  It means that 
$F(\phi) < D_{\mathcal{M}^0_0}$,\,\ $H(\phi) > 0$ and $G(\phi) < 0$ are true for all time $t$, which implies  \(\phi \not\equiv 0\). 
By applying \eqref{eqI} and \eqref{eqQ}, the conditions \(F(\phi) < D_{\mathcal{M}^0_0}\) and \(H(\phi) > 0\) lead to
\begin{align}\label{eq26}
  F(\phi)-\frac{1}{2N(p-1)}H(\phi)=\left(\frac12-\frac4{N(p-1)}\right)&\int|\Delta\phi|^2dx+\frac\omega2\int|\phi|^2dx <D_{\mathcal{M}^0_0}.
\end{align}
Due to Eq.~\eqref{eq26}, and $\frac\omega2\int|\phi|^2dx \geqslant  0$, a little effort shows that
 $$\left(\frac12-\frac4{N(p-1)}\right) \int|\Delta\phi|^2dx < D_{\mathcal{M}^0_0}.$$
According to $\frac{N+8}{N} < p <\frac{N+4}{ (N-4)^{+} }$, yield that
$\|\Delta\phi\|_2^2 < D_{\mathcal{M}^0_0}$ is true for all time $t$. 
Thus $\int|\Delta\phi|^2dx$ is bounded for any $t \in [0, T)$, so that $\phi$ is bounded in $H^2$.

Now deal with Case 2): \(\phi_0 \in \mathcal{W}^{0}_{g_2}\).  
By the assumption \(\phi_0 \in \mathcal{W}^{0}_{g_2}\) and Proposition \ref{prop3.8}, the solution \(\phi(t,x)\) of the Cauchy problem \eqref{eq1} satisfies \(\phi(t,\cdot) \in \mathcal{W}^{0}_{g_2}\) for all \(t \in [0,T)\).  
Consequently, one has $F(\phi) < D_{\mathcal{M}^0_0}$ and $ G(\phi) > 0$.
From \(G(\phi) > 0\) and \eqref{eqS}, it follows that
\begin{equation}\label{eq27}
  -\frac{\sigma}{p+1} \int|\phi|^{p+1} \, dx 
  > -\frac{1}{p+1} \left( \int|\Delta \phi|^2 \, dx 
    + \omega \int|\phi|^2 \, dx \right).
\end{equation}
On the other hand, the inequality \(F(\phi) < D_{\mathcal{M}^0_0}\) yields
\begin{equation}\label{eq28}
  \frac12 \int|\Delta \phi|^{2} \, dx 
  + \frac{\omega}{2} \int|\phi|^{2} \, dx
  - \frac{\sigma}{p+1} \int|\phi|^{p+1} \, dx
  < D_{\mathcal{M}^0_0}.
\end{equation}
Combining \eqref{eq27} and \eqref{eq28} gives
\begin{equation}\label{eq29}
  \left( \frac12 - \frac{1}{p+1} \right)
  \left( \int|\Delta \phi|^2 \, dx 
    + \omega \int|\phi|^2 \, dx \right)
  < F(\phi) < D_{\mathcal{M}^0_0}.
\end{equation}
Since $\frac{N+8}{N} < p < \frac{N+4}{(N-4)^{+}}$,  $\phi(x,t)$ remains bounded in $H^2$ for all $t \in [0,T)$.

Through the arguments of Case\:1), Case\:2) and the expression in Lemma \ref{prop2.1}, one can conclude that the solution $\phi$ to the Cauchy problem \eqref{eq1} exists globally for all $t \in [0,\infty)$, provided that $\phi_{0} \in \mathcal{W}^{0}_{g_{1}} \cup \mathcal{W}^{0}_{g_{2}}$, $\frac{N+8}{N} < p < \frac{N+4}{(N-4)^{+}}$, and $N \geqslant 2$.

$\hfill \square$

Next, we turn our attention to the mass-critical case, namely $p=\frac{N+8}{N}$.

{\bf Proof of Corollary \ref{cor1}}: 
For the case $\phi_0\in\mathcal{W}^{0}_{g_1}$, it follows from  \eqref{eq26}  that 
\be\label{eq55}
\frac\omega2\int|\phi|^2dx <D_{\mathcal{M}^0_0}.
\ee
Put $\phi_\rho=\rho^{\frac{N^2}{2N+8}}\phi(\rho x)$.  \eqref{eqQ} and $p=\frac{N+8}{N}$ imply that
\be\label{eq56}
H(\phi_\rho)=8\rho^{\frac{16}{N+4}}\int|\Delta\phi|^2dx-\frac{16\sigma}{p+1}\int|\phi|^{p+1}dx\:.
\ee
By $H(\phi)>0$, one gets $H(\phi_{\rho})<0$ when $\rho\to0$ and $H(\phi_{\rho})>0$ when $\rho\to1$. Thus by continuity, there exists a $\rho^* \in (0, 1)$ such that $H(\phi_{\rho^{*}})=0$. According to \eqref{eqI} and \eqref{eqQ}, one has
\be\label{eq57}
F(\phi_{\rho^*})=\frac{\omega}{2}\int|\phi_{\rho^*}|^2dx=\frac{\omega}{2}(\rho^{*})^{-\frac{4N}{N+4}}\int|\phi|^2dx.
\ee
From \eqref{eq55}, it follows that
\be\label{eq58}
F\left(\phi_{\rho^*}\right)<(\rho^{*})^{-\frac{4N}{N+4}}D_{\mathcal{M}^0_0}.
\ee
Also, by  \eqref{eqS} and \eqref{eqQ}, one derives from  $H(\phi_{\rho^*})=0$ and $G(\phi)<0$ that
\begin{align}\label{eq59}
  G(\phi_{\rho^{*}})=\omega(\rho^{*})^{-\frac{4N}{N+4}}\int|\phi|^2dx-\frac{4}{N}(\rho^{*})^{\frac{16}{N+4}}\int|\Delta\phi|^2dx.
\end{align}
Hence, there are  two possibilities in the following.

Case\:1): $G( \phi_{\rho ^{* }}) < 0$;

Case\:2): $G( \phi _{\rho ^{* }}) \geqslant  0$.

Treat Case\:1) firstly. In this case, noting that $G(\phi_{\rho^*})<0$ and $H(\phi_{\rho^*})=0$, one can see $\phi_{\rho^*}\in\mathcal{M}_-^0$. Thus by \eqref{eq54} and Theorem \ref{lem12}, one has
$$F(\phi_{\rho^*})\geqslant  D_{\mathcal{M}_-^0}\geqslant  D_{\mathcal{M}^0_0}>F(\phi),$$
which implies that
$F\left(\phi\right)-F\left(\phi_{\rho^*}\right)<0$,
namely,
$$\left(1-(\rho^{*})^{\frac{16}{N+4}}\right)\int|\Delta\phi|^2dx+\left(1-(\rho^{*})^{-\frac{4N}{N+4}}\right)\omega\int|\phi|^2dx<0.$$
By \eqref{eq55}, it has
\be\label{eq60}
\int|\Delta\phi|^2dx\:<\frac{(\rho^{*})^{-\frac{4N}{N+4}}-1}{1-(\rho^{*})^{\frac{16}{N+4}}}\omega\int|\phi|^2dx\:<2\frac{(\rho^{*})^{-\frac{4N}{N+4}}-1}{1-(\rho^{*})^{\frac{16}{N+4}}}D_{\mathcal{M}^0_0},
\ee
where $\frac{(\rho^{*})^{-\frac{4N}{N+4}}-1}{1-(\rho^{*})^{\frac{16}{N+4}}}>0$ for $ 0<\rho^*<1$.

Now consider the Case\:2).   In this case,  it follows from \eqref{eq58}  that
$$F\left(\phi_{\rho^*}\right)-\frac{1}{2}G(\phi_{\rho^*})=\frac{2}{N}(\rho^{*})^{\frac{16}{N+4}}\int|\Delta\phi|^2dx<(\rho^{*})^{-\frac{4N}{N+4}}D_{\mathcal{M}^0_0}.$$
This implies that 
\be\label{eq61}
\int|\Delta\phi|^2dx<\frac{N}{2}\rho^{*-4}D_{\mathcal{M}^0_0}.
\ee
By \eqref{eq60} and \eqref{eq61},  $\int|\Delta\phi|^2dx$ remains bounded for all $t\in[0,T)$. Combining with \eqref{eq55},  shows that $\phi$ is bounded in $H^2$. Therefore, Lemma \ref{prop2.1} ensures that the solution $\phi$ to the Cauchy problem \eqref{eq1} exists globally on $[0,\infty)$ for $\phi_0\in\mathcal{W}^{0}_{g_1}$ and $ p=\frac{N+8}{N}$.

For the case $\phi_0\in\mathcal{W}^{0}_{g_2}$, it is easily obtained the conclusion by putting $p=\frac{N+8}{N}$ into estimate \eqref{eq29}. Hence, the proof of Corollary \ref{cor1} is completed.

$\hfill \square$

\vspace{0.3cm}
By Theorem~\ref{th1}-{i} and the decomposition $\mathcal{W}^{0}_{g_1} \cup \mathcal{W}^{0}_{g_2} \cup \mathcal{W}^{0}_{b} = \left\{ \phi \in H^{2} \; \big| \; F(\phi) < D_{\mathcal{M}^0_0} \right\}$,
we obtain a necessary and sufficient condition for the existence of blow-up solutions to the Cauchy problem~\eqref{eq1} in the case $\mu=0$.

\begin{thm}\label{th3}
Let $\mu=0$, $\frac{N+8}{N}\leqslant  p <\frac{N+4}{ (N-4)^{+} }$ with $p\leqslant 9$, and $N\geqslant   2$. If $\phi_0$ satisfies $F(\phi_0)<D_{\mathcal{M}^0_0} $, then the solution $\phi(t,x)$ of the Cauchy problem \eqref{eq1} blows up in finite time if and only if $\phi_0\in\mathcal{W}^{0}_{b}$.
\end{thm}
In addition, by Theorem \ref{th1}-{ii} and a scaling argument, another sufficient condition for the global existence of the solution to the Cauchy problem \eqref{eq1} can be obtained.

\begin{cor}[Small Data Criterion]\label{th4}
Let $\mu=0$, $\frac{N+8}{N} \leqslant p < \frac{N+4}{(N-4)^{+}}$, $N \geqslant 2$ and $\omega > 0$.  
Suppose $\phi_{0} \in H^{2}$ satisfies  
\begin{equation}\label{eq30}
\int |\Delta \phi_0|^{2} \, dx + \omega \int |\phi_0|^{2} \, dx < 2D_{\mathcal{M}^0_0}.
\end{equation}
Then the corresponding solution $\phi(t,x)$ to the Cauchy problem \eqref{eq1} exists globally for all $t \geqslant 0$.  
\end{cor}

\begin{proof}
 From \eqref{eqI} together with \eqref{eq30}, one immediately gets $F(\phi_0) < D_{\mathcal{M}^0_0}$.  
Moreover, $G(\phi_0) > 0$ must hold; otherwise, by \eqref{eqS}, there exists some $\lambda \in (0,1]$ with $G(\lambda \phi_0) = 0$.  
Using \eqref{eqN}, Theorem \ref{lem12}, and \eqref{eq30}, one deduces that
  \be\label{eq33}
  F\left(\lambda\phi_0\right)\geqslant  D_{\mathcal{M}^0_0}.
  \ee
On the other hand, a little effort shows that
\begin{align}\label{eq31}
  \int|\Delta(\lambda\phi_0)|^2dx+\omega\int|\lambda\phi_0|^2dx=&\lambda^2\left(\int|\Delta\phi_0|^2dx+\omega\int|\phi_0|^2dx\right)\nonumber\\
  <&2\lambda^2D_{\mathcal{M}^0_0}\leqslant 2D_{\mathcal{M}^0_0},
\end{align}
which implies that
\be\label{eq32}
F\left(\lambda\phi_0\right)<D_{\mathcal{M}^0_0}.
\ee
Now, one can reach a contradiction by comparing  Eqs.~\eqref{eq33} and \eqref{eq32}. Thus, one gets $G(\phi_0)>0$. From the preceding discussion, it follows that \(\phi_0 \in \mathcal{W}^{0}_{g_2}\).  
Therefore, Theorem~\ref{th1}-{ii} directly implies the conclusion stated in Corollary~\ref{th4}.

\end{proof}
\subsection{Sharp Criterion in the case $\mu > 0$}
This subsection introduces a mixed dispersion regime due to the presence of both second-order ($\Delta$) and fourth-order ($\Delta^2$) dispersion terms. The interaction between these dispersive mechanisms breaks the standard scaling invariance, significantly complicating the associated variational structure. Here, we provide a new variational argument to study the nonlinear Schr\"odinger equations with a mixed dispersion.

To begin, a sufficient criterion for the finite-time blow-up of solutions to the Cauchy problem \eqref{eq1} is provided.
We further define
\begin{align}
&\mathcal{W}^{+}_{b}:=\left\{\phi\in H^2,\quad F(\phi)<D_{\mathcal{M}^\mu_0},\quad H(\phi)<0\right\}\label{eqKb},\\
&\mathcal{W}^{+}_{g}:=\left\{\phi\in H^{2},\quad F(\phi)<D_{\mathcal{M}^\mu_0}, \quad H(\phi)>0\right\}\label{eqK+}.
\end{align}
and $\phi_{\lambda}(x)=\lambda^{\frac{N}{2}}\phi(\lambda x)$. Thus the following necessary Lemmas hold.
\begin{lem}\label{lemdi1}
  Let $\mu>0$ and $\phi\in H^2$. When $p = \frac{N+8}{N}$, we further assume that $H(\phi)\leqslant0$. Then there exists $\lambda_0>0$ such that $G(\phi_{\lambda_0}) = 0$.
  \end{lem}
  \begin{proof}
  Noting that $G(\phi_{\lambda})$ is expressed as
  \[
    G(\phi_{\lambda})=\lambda^4\int|\Delta \phi|^2dx+\lambda^2\mu\int|\nabla \phi|^2dx+\omega\int|\phi|^2dx-\lambda^{\frac{N(p-1)}{2}}\int|\phi|^{p + 1}dx.
  \]
  If $p>\frac{N+8}{N}$, that is, $\frac{N(p-1)}{2}>4$, it can be seen that this result is a direct consequence of the structure of the graph of the mapping \(\lambda \mapsto G(\phi_{\lambda})\).
  
  If $p=\frac{N+8}{N}$, one has the expression
\be\label{eq81}
    G(\phi_{\lambda})=\lambda^4(\int|\Delta \phi|^2dx-\sigma\int|\phi|^{p + 1}dx)+\lambda^2\mu\int|\nabla \phi|^2dx+\omega\int|\phi|^2dx.
\ee
  From $H(\phi)\leqslant0$, it follows that
 $$ \int|\Delta \phi|^{2}dx+\frac{\mu}{2}\int|\nabla \phi|^{2}dx\leqslant \frac{2\sigma}{p+1}\int|\phi|^{p+1}dx.$$
Therefore, one can deduce to
  \[
    \int|\Delta \phi|^2dx-\sigma\int|\phi|^{p + 1}dx<\int|\Delta \phi|^2dx-\frac{2\sigma}{p+1}\int|\phi|^{p + 1}dx\leqslant-\frac{\mu}{2}\int|\nabla \phi|^{2}dx<0.
  \]
So the highest power coefficient of \(\lambda\) in Eq.~\eqref{eq81} is negative. This completes the proof.

  \end{proof}
  \begin{lem}\label{lemdi2}
    Assume that $\mu>0$, $\phi\in H^2$ satisfies $\phi\neq0$ and $H(\phi)\leqslant0$. Then
    \[
    \frac{1}{16}H(\phi)\leqslant F(\phi)-D_{\mathcal{M}^\mu_0}.
    \]
    \end{lem}
    \begin{proof}
    Define the function
    \begin{align*}
    f(\lambda):&=F(\phi_{\lambda})-\frac{\lambda^4}{16}H(\phi)\\
    &=\frac{\omega}{2}\|\phi\|_{L^2}^2+\frac{\mu\lambda^2}{2}\|\nabla \phi\|_{L^2}^2 +\frac{1}{2}\lambda^4(\|\Delta \phi\|_{L^2}^2-\frac{1}{8}H(\phi))-\frac{\sigma\lambda^{\frac{N(p-1)}{2}}}{p + 1}\|\phi\|_{L^{p + 1}}^{p + 1}.
    \end{align*}
    First, one proves the important process $f(1)=\max\limits_{\lambda>0}f(\lambda)$. To this aim,
    the analysis are  divided into two cases.
    
   Case 1): When $p>\frac{N+8}{N}$. Since $\frac{N(p-1)}{2}>4$, after some  calculations,  one deduces that  
   \begin{align}
    \frac{\partial f(\lambda)}{\partial \lambda}&=2\lambda^3(\|\Delta \phi\|_{L^2}^2-\frac{1}{8}H(\phi))+\mu\lambda\|\nabla \phi\|_{L^2}^2 -\frac{N\sigma(p-1)}{2(p+1)}\lambda^{\frac{N(p-1)}{2}-1}\|\phi\|_{L^{p + 1}}^{p + 1}\nonumber\\
&=\lambda\big[\lambda^2(\frac{N\sigma(p-1)}{2(p+1)}\|\phi\|_{L^{p + 1}}^{p + 1}-\mu\|\nabla \phi\|_{L^2}^2)+\mu\|\nabla \phi\|_{L^2}^2 -\lambda^{\frac{N(p-1)}{2}-2}\frac{N\sigma(p-1)}{2(p+1)}\|\phi\|_{L^{p + 1}}^{p + 1}\big]\nonumber\\
 &=\lambda\big[\lambda^2(1-\lambda^{\frac{N(p-1)}{2}-4})\frac{N\sigma(p-1)}{2(p+1)}\|\phi\|_{L^{p + 1}}^{p + 1}+(1-\lambda^2)\mu\|\nabla \phi\|_{L^2}^2\big],
   \end{align}
which implies that  $f^\prime(1) = 0$. Hence, the graphical analysis of \(f(\lambda)\) reveals that $f(1)=\max\limits_{\lambda>0}f(\lambda)$.
   
   \vspace{0.2cm}
  Case 2): When $p = \frac{N+8}{N}$, the expression of $f(\lambda)$ can be written as    \begin{align*}
      f(\lambda)&=\frac{\omega}{2}\|\phi\|_{L^2}^2+\frac{\mu\lambda^2}{2}\|\nabla \phi\|_{L^2}^2 +\frac{1}{2}\lambda^4[\|\Delta \phi\|_{L^2}^2-\frac{1}{8}H(\phi)-\frac{2\sigma}{p + 1}\|\phi\|_{L^{p + 1}}^{p + 1}]\\
      &=\frac{\omega}{2}\|\phi\|_{L^2}^2+\frac{\mu\lambda^2}{2}\|\nabla \phi\|_{L^2}^2 +\frac{1}{2}\lambda^4[-\frac{\mu}{2}\|\nabla \phi\|_{L^2}^2].
    \end{align*}
Since \(f'(1) = 0\) and considering the shape of the graph of \(f(\lambda)\), it follows that $f(1)=\max\limits_{\lambda>0}f(\lambda)$.
    
  Next, the desired inequality will be established.  
By Lemma~\ref{lemdi1}, there exists \(\lambda_0 > 0\) such that \(G(\phi_{\lambda_0}) = 0\). Hence, by the definition of \(D_{\mathcal{M}^\mu_0}\), $H(\phi)\leqslant0$, and $f(1)\geqslant f(\lambda_0)$, one obtains
    \[
      D_{\mathcal{M}^\mu_0}\leqslant F(\phi_{\lambda_0})\leqslant F(\phi_{\lambda_0})-\frac{\lambda_0^4}{16}H(\phi)\leqslant F(\phi)-\frac{1}{16}H(\phi).
    \]

    Hence, this completes the proof.

    \end{proof}
   \begin{lem}\label{prop3.7}
Let $\mu>0$, $\frac{N+8}{N} \leqslant p < \frac{N+4}{(N-4)^{+}}$ and $N \geqslant 2$.  
Then the sets $\mathcal{W}^{+}_{b}$ and $\mathcal{W}^{+}_{g}$ are invariant manifolds for the Cauchy problem \eqref{eq1}.  
More precisely, if $\phi_0 \in \mathcal{W}^{+}_{b}$, then the corresponding solution $\phi(t,x)$ of \eqref{eq1} satisfies
$\phi(t, \cdot) \in \mathcal{W}^{+}_{b}$ for all $t \in [0, T)$.
Similarly, the same invariance property holds for $\mathcal{W}^{+}_{g}$.
\end{lem}

      \begin{proof} 
    Let \(\phi_0 \in \mathcal{W}^{+}_{b}\).  
By Lemma~\ref{prop2.1}, there exists a unique solution \(\phi(t, \cdot) \in C([0,T); H^2)\) with \(T \leqslant \infty\) to the Cauchy problem \eqref{eq1} such that \(\phi(t,x)\) satisfies the Eq.~\eqref{eq1}.  
From the conservation laws \eqref{eqm} and \eqref{eqe1}, one sees that 
\begin{equation}\label{eq17}
F(\phi(t, \cdot)) = F(\phi_0), \quad \text{for all } t \in [0,T).
\end{equation}
In particular, the condition \(F(\phi_0) < D_{\mathcal{M}^\mu_0}\) guarantees  $F(\phi(t, \cdot)) < D_{\mathcal{M}^\mu_0}$ for all $ t \in [0, T)$.     

    Next, the inequality \(H(\phi(t,\cdot)) < 0\) for all \(t \in [0,T)\) is established by contradiction.  
Assume, on the contrary, that there exists some \(t^{*} \in (0,T)\) such that  
\begin{equation*}
    H(\phi(t^{*}, \cdot)) = 0 \quad and \quad
    0 \leqslant F(\phi(t^{*})) - D_{\mathcal{M}^\mu_0},
\end{equation*}
where the latter follows from Lemma~\ref{lemdi2}. This contradicts the conservation law 
$F(\phi(t^{*})) = F(\phi_0) < D_{\mathcal{M}^\mu_0}$.
Therefore, it must hold that \(H(\phi(t,\cdot)) < 0\) for all \(t \in [0,T)\).

From the above arguments, it follows that \(\phi(t,\cdot) \in \mathcal{W}^{+}_{b}\) for all \(t \in [0,T)\).  
By an analogous argument, the invariance of \(\mathcal{W}^{+}_{g}\) can also be established.

      \end{proof}
Relying on the above lemmas, we can prove Theorem \ref{thmdi1} in this section systematically.



\vspace{0.3cm}
{\bf Proof of Theorem \ref{thmdi1}-{i}:}
According to \(\phi_0\in \mathcal{W}^{+}_{b}\) and define $$b:=D_{\mathcal{M}^\mu_0}-F(\phi_0)>0.$$ Similar to the proof of Theorem \ref{th1}-{i}, the key to this proof lies in obtaining the inequality
\be\label{eq73.2}
\frac{d}{dt} \mathcal{M}_{\varphi_R}[\phi(t)] \leqslant -\tau \|\nabla \phi(t)\|^{2}_{2}<0, \quad \text{for all } \,\ t\in [0, T).
\ee
To achieve this, one discusses it in the following two cases.

Case 1): Assume that
$$\|\nabla \phi(t)\|^2_2 \leqslant \frac{4N(p-1) E(\phi_0)}{\mu(N(p-1) - 4)}$$
holds for all time  $ t \in [0, T)$. By choosing \(R > 0\) sufficiently large, one can require the following estimate,
$$O\left( R^{-4} + \|\nabla\phi(t)\|^2_2 R^{-2} + \|\nabla\phi(t)\|^{\frac{p-1}{2}}_2 R^{-\frac{(p-1)(N-1)}{2}} + |\mu| R^{-2} \right) < b.$$
From Lemma \ref{prop2.1}, Lemma \ref{lem11} and Lemma \ref{lemdi2}, a little careful calculation  shows that
    \begin{align*}
    \frac d{dt}\mathcal{M}_{\varphi_R}[\phi(t)]&\leqslant H(\phi)+O\left( R^{-4} + \|\nabla\phi(t)\|^2_2 R^{-2} + \|\nabla\phi(t)\|^{\frac{p-1}{2}}_2 R^{-\frac{(p-1)(N-1)}{2}} + |\mu| R^{-2} \right) \\
    &\leqslant 16(F(\phi_0)-D_{\mathcal{M}^\mu_0})+b\\
    &=-15b.
    \end{align*}
Thus,
\begin{equation}\label{eq71.1}
\frac{d}{dt}\mathcal{M}_{\varphi_R}[\phi(t)] \leqslant -15b \leqslant -\tau\|\nabla \phi(t)\|^2_2
\end{equation}
provided $\tau > 0$ is sufficiently small.

Case 2): Suppose that  for all $ t \in [0, T)$, 
$$\|\nabla \phi(t)\|^2_2 > \frac{4N(p-1) E(\phi_0)}{\mu(N(p-1) - 4)}$$
holds. According to Lemma \ref{lem11} and $$\|\nabla\phi(t)\|^2_2 \leqslant \|\phi(t)\|_2\|\Delta\phi(t)\|_2 =\|\phi(0)\|_2\|\Delta\phi(t)\|_2,$$
 one derives
\begin{align*}
\frac{d}{dt} \mathcal{M}_{\varphi_R}[\phi(t)] & \leqslant 2N(p-1) E(\phi_0) - (N(p-1) - 8)\|\Delta \phi(t)\|_{2}^2-(N(p-1) - 4)\mu\|\nabla\phi(t)\|_{2}^2\\
&\quad +O\left(R^{-4} + R^{-2}\|\nabla \phi(t)\|_{2}^2 + R^{-\frac{(p-1)(N-1)}{2}}\|\nabla \phi(t)\|_{2}^{\frac{p-1}{2}} + |\mu|R^{-2}\right)\\
&\leqslant  -\frac{1}{2}(N(p-1) - 4)\mu\|\nabla\phi(t)\|_{2}^2 - (N(p-1) - 8)\frac{\|\nabla \phi(t)\|_{2}^4}{\|\phi(0)\|_{2}^2}\\
&\quad +O\left(R^{-4} + R^{-2}\|\nabla \phi(t)\|_{2}^2 + R^{-\frac{(p-1)(N-1)}{2}}\|\nabla \phi(t)\|_{2}^{\frac{p-1}{2}} + |\mu|R^{-2}\right).
\end{align*}
By taking \(R > 0\) sufficiently large and $p\leqslant 9$, it follows that
\begin{equation}\label{eq72.1}
\frac{d}{dt}\mathcal{M}_{\varphi_R}[\phi(t)] \leqslant -\frac{1}{4}(N(p-1) - 4)\mu\|\nabla\phi(t)\|_{2}^2.
\end{equation}
Combining the estimates \eqref{eq71.1} and \eqref{eq72.1}, namely exists a $\tau>0$ such that
$\frac{d}{dt} \mathcal{M}_{\varphi_R}[\phi(t)] \leqslant -\tau \|\nabla \phi(t)\|^{2}_{2}<0$ for all $t\in [0, T)$.

By the monotonicity of the function $\mathcal{M}_{\varphi_R}[\phi]$, there exists $t_2 \geqslant 0$ such that 
\be\label{eq74}
\mathcal{M}_{\varphi_R}[\phi(t)] < 0, \quad \text{for} \,\ t \geqslant t_2.
\ee

Next, it can be shown that the solution $\phi ( t, x)$ to the Cauchy problem \eqref{eq1} blows up in a  finite time.
Integrating \eqref{eq73.2} on $[t_2, t]$ and using \eqref{eq74}, one gets
$$\mathcal{M}_{\varphi_R}[\phi(t)] \leqslant  -\tau \int_{t_2}^t\|\nabla \phi(s)\|^{2}_{2}ds.$$
Injecting the Cauchy--Schwarz inequality into $\mathcal{M}_{\varphi_R}[\phi(t)]$ yields   
$$
|\mathcal{M}_{\varphi_R}[\phi(t)]| \leqslant 2\|\nabla\varphi_R\|_\infty\|\phi(t)\|_2\|\nabla\phi(t)\|_2 \leqslant C(\phi_0)R\|\nabla\phi(t)\|_2,
$$
and therefore
$$\int_{t_2}^t |\mathcal{M}_{\varphi_R}[\phi(s)]|^2ds \leqslant C^2(\phi_0)R^2\int_{t_2}^t\|\nabla\phi(s)\|_2^2ds \leqslant -\frac{C^2(\phi_0)R^2}{\tau}\mathcal{M}_{\varphi_R}[\phi(t)].$$
Thus, denoting $B:=\frac{\tau}{C^2(\phi_0)R^2}$ and
 $\gamma_2 (t):=\int_{t_2}^t |\mathcal{M}_{\varphi_R}[\phi(s)]|^2ds$ for $t\in[t_2,T)$, one has
$\mathcal{M}_{\varphi_R}[\phi(t)]\leqslant -B\gamma_2(t)$.  The function $\gamma_2$ is
nonnegative and strictly increasing by \eqref{eq74}. Moreover,
$\gamma_2'(t)=\big|\mathcal{M}_{\varphi_R}[\phi(t)]\big|^2\geqslant B^2\gamma_2(t)^2$. Fixing a
time $t_3\in(t_2,T)$, so that $\gamma_2(t_3)>0$, and integrating this differential inequality
on $[t_3,t]$, one can deduce that
\[
\gamma_2 (t) \geqslant \frac{\gamma_2 (t_3)}{1 - B^2(t - t_3)\gamma_2(t_3)},
\]
and thus $\mathcal{M}_{\varphi_R}[\phi(t)] \to -\infty$ as $t$ tends to $t^*_2:=t_3+\frac{1}{B^2 \gamma_2(t_3)}$. Since
$|\mathcal{M}_{\varphi_R}[\phi(t)]|\leqslant C(\phi_0)R\|\nabla\phi(t)\|_2$ and
$\|\nabla\phi(t)\|_2^2\leqslant\|\phi_0\|_2\|\Delta\phi(t)\|_2$, this forces
$\|\Delta\phi(t)\|_2\to+\infty$, and hence $T\leqslant t^*_2<+\infty$, which completes the proof under the assumptions \(\mu > 0\) and \(\frac{N+8}{N} \leqslant p < \frac{N+4}{(N-4)^{+}}\).

$\hfill \square$

\vspace{0.3cm}
{\bf Proof of Theorem \ref{thmdi1}-{ii}:} For \(\phi_0\in \mathcal{W}^{+}_{g}\) and \(\phi(t,x)\) be the solution of \eqref{eq1} corresponding the initial datum \(\phi_0\) in \(t\in[0,T)\).
    It follows from Lemma \ref{prop3.7} that \(\phi\in \mathcal{W}^{+}_{g}\), that is
    \begin{equation}
      F(\phi)=\frac{1}{2}T_{\omega,\mu}(\phi)-\frac{\sigma}{p+1}\int|\phi|^{p+1}dx <D_{\mathcal{M}^\mu_0},\label{eqdi9}
    \end{equation}
    and 
    \begin{equation}
    H(\phi)=8\int|\Delta \phi|^{2}dx+4\mu\int|\nabla \phi|^{2}dx-\frac{2N\sigma(p-1)}{p+1}\int|\phi|^{p+1}dx>0\label{eqdi10}.
  \end{equation}
    By \eqref{eqdi9} and \eqref{eqdi10}, one deduces
    \[
      D_{\mathcal{M}^\mu_0}>F(\phi)>\int\frac{N(p - 1)-8}{2N(p - 1)}|\Delta\phi|^2+\frac{\mu (N(p - 1)-4)}{2N(p - 1)}|\nabla \phi|^2+\frac{\omega}{2}|\phi|^2dx,
    \]
    which yields $\frac{N+8}{N}< p<\frac{N+4}{ (N-4)^{+} }$ when \(N\geqslant2\). Therefore, from Lemma \ref{prop2.1}, the solution \(\phi(t,x)\) of the Cauchy problem \eqref{eq1} globally exists on \(t\in[0,\infty)\).

    $\hfill \square$

\begin{rem}
  It is worth emphasizing that the approach used for the case $\mu = 0$ cannot be extended to the case $\mu > 0$ directly. 
   For $\mu = 0$, we adopt two-constraint variational framework based on the sets 
$\mathcal{M}^\mu_0$
  and 
$\mathcal{M}_-^0$
 , which provides sharp threshold for blow-up. However, when 
$\mu > 0$, the two-constraint method becomes substantially more difficult to apply, mainly due to the fact that the constraint set 
$\mathcal{M}_-^\mu$ may be empty, that is, it becomes extremely challenging to find a scaling such that both
$G(u)<0$ and $H(u)=0$ hold simultaneously.

To overcome this, we instead employ one-constraint variational approach based only on 
$D_{\mathcal{M}^\mu_0}$, combined with method inspired by Hayashi et al \cite{Hayashi}, to derive sharp threshold condition for 
$\mu > 0$. While this one-constraint method still yields valid and sharp threshold for blow-up and global existence, the comparison between the two methods clearly indicates that the two-constraint method produces more accurate results when $\mu = 0$.

In conclusion, the sign of $\mu$ determines two qualitatively distinct analytical frameworks. Each requires specific functional analytic and variational approaches to handle the corresponding sharp threshold problem.
\end{rem}

\newpage

\renewcommand{\theequation}{\thesection.\arabic{equation}}
\setcounter{equation}{0}
\section{Instability of the Standing Waves}
By Lemma~\ref{lem14}, the variational problem
\be\label{eq42}
D_{\mathcal{M}^\mu_0}=\min_{u\in\mathcal{M}^\mu_0}F\left(u\right)
\ee
admits a minimizer $u \in \mathcal{M}^\mu_0$. 
The attainment of this minimum implies that $u$ satisfies the corresponding Euler--Lagrange equation, namely the elliptic equation \eqref{eq64}, and hence, as a minimizer, $u$ is a ground state solution of \eqref{eq64}.
Consequently, $\phi(t,x) = e^{i\omega t} u(x)$ is a standing wave solution of \eqref{eq1}. 
Now, prior to proving the instability result, several auxiliary lemmas are established, which serve as key tools in the subsequent analysis.
\begin{lem}\label{lem8}
Suppose that $\phi \in H^2 \setminus \{0\}$ and $\mu \geqslant 0$. There exists a unique positive constant $\kappa$ such that $G(\kappa \phi) = 0$
and moreover, $F(\kappa \phi)$ attains the global maximum along the ray $\{\lambda \phi : \lambda > 0\}$.
\end{lem}
\begin{proof} For $\lambda>0$. It following from the scaling that 
\begin{align}
  G(\lambda\phi)&=\lambda^{2}T_{\omega,\mu}(\phi)-\sigma\lambda^{p+1}\int|\phi|^{p+1}dx\:,\label{eq45}\\
  F\left(\lambda\phi\right)&=\frac{1}{2}\lambda^{2}T_{\omega,\mu}(\phi)-\frac{1}{p+1}\lambda^{p+1}\sigma\int|\phi|^{p+1}dx\label{eq46}.
\end{align}
Noting that \eqref{eq45} and \eqref{eq46}, then
\begin{align}\label{eq47}
\frac{d}{d\lambda}F\left(\lambda\phi\right)&=\lambda T_{\omega,\mu}(\phi)-\sigma\lambda^{p}\int|\phi|^{p+1}dx =\frac{G(\lambda\phi)}{\lambda}.
\end{align}
It can be known from \eqref{eq45} that $G(\lambda\phi)>0$ when $\lambda\to 0$, and $G(\lambda\phi)\leqslant0$ when $\lambda\to +\infty$. Then, there  exists a unique $\kappa>0$ such that $G(\kappa \phi)=0$, and Eq.~\eqref{eq47} implies $F'(\kappa \phi) = 0$. Consider  the qualitative behavior of the function $\lambda \mapsto F(\lambda \phi)$. It is clear that
\[
F(\kappa \phi) = \max_{\lambda > 0} F(\lambda \phi).
\]  
Hence, the statement of Lemma~\ref{lem8} is verified.

\end{proof}

\begin{lem}\label{lem9}
Let $u$ be a minimizer of \eqref{eqN}. Then $H(u)=0$. 
\end{lem}
\begin{proof}
As $u$ attains the minimum in \eqref{eqN}, it also satisfies the elliptic equation \eqref{eq64}.  
Multiplying \eqref{eq64} by $x \cdot \nabla u$ yields
\be\label{eq48}
\left(2 - \frac{N}{2}\right)\int |\Delta u|^2dx + \mu\left(1 - \frac{N}{2}\right)\int |\nabla u|^2dx - \frac{\omega N}{2}\int |u|^2dx + \frac{\sigma N}{p+1}\int |u|^{p+1}dx = 0.
\ee
which is the Pohozaev identity of Eq.~\eqref{eq64}. Note that $G(u)=0$ and \eqref{eq48} implies that $H( u) = 0$. 

\end{proof}
Applying Lemmas \ref{lem8} and \ref{lem9} to the standing wave $\phi(t,x) = e^{i\omega t} u(x)$ leads to the following theorem on the instability of Eq.~\eqref{eq1}.

\vspace{0.3cm}
{\bf Proof of Theorem \ref{th5.1}}:    First of all, 
under  the assumption that  $p$ is odd integer satisfying $ \frac{N+8}{N}\leqslant p < \min\left\{\frac{N+4}{(N-4)^{+}}, 9\right\} $,  one can 
apply the Fourier rearrangement (see also Appendix A. 2 in  \cite{Boulenger2017}) to deduce that the minimizer on  $M_0^{\mu}$  can be taken to be radial.
Then, for $\mu=0$, one can take  $u \in H^2 \setminus \{0\}$ be a radial  minimizer of \eqref{eqN} as ensured by Lemma \ref{lem14} and Lemma \ref{lem3.2}, satisfying $G(u) = 0$ and $H(u)=0$.
It follows from \eqref{eqS} and \eqref{eqQ} that
\begin{align}
  G(\lambda u)&=\lambda^{2}\left(\int|\Delta u|^{2}dx+\omega\int|u|^{2}dx-\sigma\lambda^{p-1}\int|u|^{p+1}dx\right)\:,\\
  H(\lambda u)&=\lambda^{2}\left(8\int|\Delta u|^{2}dx-\frac{2N\sigma(p-1)}{p+1}\lambda^{p-1}\int|u|^{p+1}dx\right).
\end{align}
One can deduce that,
\be\label{eq49.1}
G(\lambda u)<0,\quad H(\lambda u)<0,\quad  \text{for any} \,\ \lambda>1.
\ee
On the other hand, combining \eqref{eq49.1} with $G(u)=0$ and Lemma \ref{lem8} show that $\kappa=1$ and $F(\lambda u)<F(u)$ whenever $\lambda>1$. Since $F(u)=D_{\mathcal{M}^0_0}$, this implies $\lambda u \in \mathcal{W}^{0}_{b}$ for any $\lambda>1$.

Next, let  the initial data $\phi_0 = \lambda u(x)$  
with $\lambda$ chosen sufficiently close to 1 so that
$$\|\phi_0-u\|_{H^2}=\|\lambda u-u\|_{H^2}=(\lambda -1)\|u\|_{H^2}<\varepsilon. $$
By Theorem \ref{th1}-{i}, the solution $\phi(t,x)$ of the Cauchy problem \eqref{eq1} blows up in a finite time. Therefore, according to the definition of instability for standing waves, the corresponding solution $u \in H^2$ of \eqref{eq64} is unstable.

Similar to $\mu=0$, one can obtain the proof of $\mu>0$.    For $\mu>0$, let $\tilde{u} \in H^2 \setminus \{0\}$ be a radial minimizer of \eqref{eqN} satisfying $G(\tilde{u}) = 0$ and $H(\tilde{u})=0$. Furthermore, for any $\lambda > 1$, it holds that
$G(\lambda \tilde{u})<0$, $H(\lambda \tilde{u})<0$.
On the basis of Lemma \ref{lem8}, it means that $F(\lambda \tilde{u})<F(\tilde{u})=D_{\mathcal{M}^\mu_0}$. Thus one can obtain that $\lambda \tilde{u}\in\mathcal{W}^{+}_{b}$ for any $\lambda>1$. 

Then let $\tilde{\phi}_0=\lambda \tilde{u}(x)$, where $\lambda>0$ is taken near 1 such that
$$\|\phi_0-\tilde{u}\|_{H^2}=\|\lambda \tilde{u}-\tilde{u}\|_{H^2}=(\lambda -1)\|\tilde{u}\|_{H^2}<\varepsilon. $$
By Theorem \ref{thmdi1}-{i}, the solution $\phi(t,x)$ of the Cauchy problem \eqref{eq1} blows up in a finite time. Consequently, the corresponding standing wave $\tilde{u} \in H^2$ of \eqref{eq64} is unstable.

In conclusion, the statement is proved.

$\hfill \square$


\section*{Declarations}

\qquad {\bf Funding}

This work was supported by  
the National Natural Science Foundation of China
(Grant No. 12071323).

 \vspace{0.4cm}

{\bf Data Availability Statement}

Data will be made available on request.

 \vspace{0.4cm}
{\bf Conflict of Interest}

The author declares that there are no competing inter-
ests.
 
 \vspace{0.4cm}
{\bf  Ethical approval}

This article contains no studies involving human participants or animals performed by the authors. Ethical approval was therefore not
required.

 \vspace{0.4cm}
{\bf Author contributions}

All authors Contributed to the study conception and design. Material preparation, data collection and analysis were performed by Shanshan Fu and Shihui Zhu.
The first draft of the manuscript was written by Shanshan Fu and Shihui Zhu. All authors
commented on previous versions of the manuscript, and all authors read and approved the final
manuscript.

 \vspace{0.4cm}
{\bf Data availability}

 No datasets were generated or analyzed in this study. Data
sharing is therefore not applicable.

 \vspace{0.4cm}
{\bf Acknowledgements}  

The author sincerely thanks the handling editor and
the anonymous referees for their careful evaluation and constructive comments,
which helped improve the clarity and quality of the manuscript.

 \vspace{1.5cm}

\noindent {\bf References}
\newcounter{cankao}
\begin{list}
{[\arabic{cankao}]}{\usecounter{cankao}\itemsep=0cm}

\bibitem{BFM2010} Baruch, G., Fibich, G., Mandelbaum, E.: Singular solutions of the biharmonic nonlinear Schr\"odinger equation. SIAM J. Appl. Math. 70(8), 3319-3341 (2010)

\bibitem{BenArtzi2000} Ben-Artzi, M., Koch, H., Saut, J.C.: Dispersion estimates for fourth order Schr\"odinger equations. C. R. Acad. Sci. Paris Sér. I Math. 330(2), 87-92 (2000)

\bibitem{BCGJ2019TAMS} Bonheure, D., Casteras, J.B., Gou, T., Jeanjean, L.: Normalized solutions to the mixed dispersion nonlinear Schr\"odinger equation in the mass critical and supercritical regime. Trans. Amer. Math. Soc. 372(3), 2167-2212 (2019)

\bibitem{BCGJ2019} Bonheure, D., Castéras, J.B., Gou, T., Jeanjean, L.: Strong instability of ground states to a fourth order Schr\"odinger equation. Int. Math. Res. Not. 2019(17), 5299-5315 (2019)

\bibitem{Boulenger2017} Boulenger, T., Lenzmann, E.: blow-up for biharmonic NLS. Ann. Sci. Éc. Norm. Supér. 50(3), 503-544 (2017)

\bibitem{CL1982} Cazenave, T., Lions, P.L.: Orbital stability of standing waves for some nonlinear Schr\"odinger equations. Comm. Math. Phys. 85, 549-561 (1982)

\bibitem{CKSTT2002} Colliander, J., Keel, M., Staffilani, G., Takaoka, H., Tao, T.: Almost conservation laws and global rough solutions to a nonlinear Schr\"odinger equation. Math. Res. Lett. 9, 659-682 (2002)

\bibitem{CuiGuo2007} Cui, S., Guo, C.: Well-posedness of higher-order nonlinear Schr\"odinger equations in Sobolev spaces \(H^{s}(\mathbb{R}^{N})\) and applications. Nonlinear Anal. 67(3), 687-707 (2007)

\bibitem{Dinh2019} Dinh, V.D.: On blowup solutions to the focusing intercritical nonlinear fourth-order Schr\"odinger equation. J. Dynam. Differential Equations 31(4), 1793-1823 (2019)

\bibitem{Dinh2021} Dinh, V.D.: Dynamics of radial solutions for the focusing fourth-order nonlinear Schr\"odinger equations. Nonlinearity 34(2), 776-821 (2021)

\bibitem{FJMM2022} Fernandez, A.J., Jeanjean, L., Mandel, R., Maris, M.: Non-homogeneous Gagliardo-Nirenberg inequalities in \(\mathbb{R}^n\) and application to a biharmonic nonlinear Schr\"odinger equation. J. Differ. Equ. 330, 1-65 (2022)

\bibitem{Fibich2002} Fibich, G., Ilan, B., Papanicolaou, G.: Self-focusing with fourth-order dispersion. SIAM J. Appl. Math. 62(4), 1437-1462 (2002)

\bibitem{Fibich2003} Fibich, G., Ilan, B., Schochet, S.: Critical exponents and collapse of nonlinear Schr\"odinger equations with anisotropic fourth-order dispersion. Nonlinearity. 16, 1809-1821 (2003)

\bibitem{Fibich2015} Fibich, G.: The nonlinear Schr\"odinger equation: singular solutions and optical collapse. Appl. Math. Sci., vol. 192. Springer, New York (2015)

\bibitem{Hayashi} Fukaya, N., Hayashi, M.: Instability of algebraic standing waves for nonlinear Schr\"odinger equations with double power nonlinearities. Trans. Amer. Math. Soc. 374(2), 1421-1447 (2021)

\bibitem{GSS1987} Grillakis, M., Shatah, J., Strauss, W.: Stability theory of solitary waves in the presence of symmetry. I. J. Funct. Anal. 74(1), 160-197 (1987)

\bibitem{HN2025} Hayashi, N., Naumkin, P.I.: Global existence of odd solutions for the cubic fourth-order nonlinear Schr\"odinger equations. NoDEA Nonlinear Differ. Equ. Appl. 32(3), 51 (2025)

\bibitem{KM2024} Komada, K., Masaki, S.: Scattering of solutions with group invariance for the fourth-order nonlinear Schr\"odinger equation. Nonlinearity. 37(8), 085003 (2024)

\bibitem{HolmerR2010} Holmer, J., Platte, R., Roudenko, S.: blow-up criteria for the 3D cubic nonlinear Schr\"odinger equation. Nonlinearity. 23(4), 977-1030 (2010)

\bibitem{KS2000} Kappan, V.I., Shagalov, G.: Stability of solitons described by nonlinear Schr\"odinger-type equations with higher-order dispersion. Phys. D 144, 194-210 (2000)

\bibitem{KenigMerle2006IM} Kenig, C.E., Merle, F.: Global well-posedness, scattering and blow-up for the energy-critical, focusing, nonlinear Schr\"odinger equation in the radial case. Invent. Math. 166, 645-675 (2006)

\bibitem{kenig1991} Kenig, C.E., Ponce, G., Vega, L.: Oscillatory integrals and regularity of dispersive equations. Indiana Univ. Math. J. 40(1), 33-69 (1991)

\bibitem{K2022} Komada, K.: Existence of blow-up solutions to nonlinear Schr\"odinger equations with anisotropic fourth-order dispersion. Funkcial. Ekvac. 65(2), 191-214 (2022)

\bibitem{LW2024} Lenzmann, E., Weth, T.: Symmetry breaking for ground states of biharmonic NLS via Fourier extension estimates. J. Anal. Math. 152(2), 777-800 (2024)

\bibitem{Lin2022} Lin, Z., Zeng, C.: Instability, index theorem, and exponential trichotomy for linear Hamiltonian PDEs. Mem. Amer. Math. Soc. 275, 1347 (2022)

\bibitem{LuoZhengZhu2023} Luo, T., Zheng, S., Zhu, S.: The existence and stability of normalized solutions for a bi-harmonic nonlinear Schr\"odinger equation with mixed dispersion. Acta Math. Sci. 43B(2), 539-563 (2023)

\bibitem{Ma2023} Ma, L.: Positive radially symmetric ground states to a bi-harmonic problem. Differ. Integral Equ. 36(3/4), 313-320 (2023)

\bibitem{Merle2005} Merle, F., Raphael, P.: The blow-up dynamic and upper bound on the blow-up rate for critical nonlinear Schr\"odinger equation. Ann. Math. 161, 157-222 (2005)

\bibitem{MZ2016} Miao, C.X., Zheng, J.Q.: Scattering theory for the defocusing fourth-order Schr\"odinger equation. Nonlinearity. 29(2), 692-736 (2016)

\bibitem{NP2015} Natali, F., Pastor, A.: The fourth-order dispersive nonlinear Schr\"odinger equation: orbital stability of a standing wave. SIAM J. Appl. Dyn. Syst. 14(3), 1326-1347 (2015)

\bibitem{Pausader2007} Pausader, B.: Global well-posedness for energy critical fourth-order Schr\"odinger equations in the radial case. Dyn. Partial Differ. Equ. 4(3), 197-225 (2007)

\bibitem{Pausader2007b} Pausader, B.: Scattering and the Levandosky-Strauss conjecture for fourth-order nonlinear wave equations. J. Differ. Equ. 241(2), 237-278 (2007)

\bibitem{PX2013} Pausader, B., Xia, S.: Scattering theory for the fourth-order Schr\"odinger equation in low dimensions. Nonlinearity. 26(8), 2175-2191 (2013)

\bibitem{SS2020} Saut, J.C., Segata, J.: Long range scattering for the nonlinear Schr\"odinger equation with higher order anisotropic dispersion in two dimensions. J. Math. Anal. Appl. 483(2), 123638 (2020)

\bibitem{Weinstein1983} Weinstein, M.I.: Nonlinear Schr\"odinger equations and sharp interpolation estimates. Commun. Math. Phys. 87, 567-576 (1983)

\bibitem{YYZ2023} Yu, X.Y., Yue, H.T., Zhao, Z.H.: On the decay property of the cubic fourth-order Schr\"odinger equation. Proc. Amer. Math. Soc. 151(6), 2619-2630 (2023)

\bibitem{YYZ2024} Yu, X.Y., Yue, H.T., Zhao, Z.H.: Global well-posedness and scattering for fourth-order Schr\"odinger equations on waveguide manifolds. SIAM J. Math. Anal. 56(1), 1427-1458 (2024)

\bibitem{Zhang2002} Zhang, J.: Sharp conditions of global existence for nonlinear Schr\"odinger and Klein-Gordon equations. Nonlinear Anal. 48(2), 191-207 (2002)

\bibitem{Zhang2005} Zhang, J.: Sharp threshold for blow-up and global existence in nonlinear Schr\"odinger equations under a harmonic potential. Comm. Partial Differ. Equ. 30, 1429-1443 (2005)

\bibitem{ZhangWang2024} Zhang, Z., Wang, D., Chen, J., Xie, Z., Xu, C.: blow-up of solutions for the fourth-order Schr\"odinger equation with combined power-type nonlinearities. J. Geom. Anal. 34(10), 307 (2024)

\bibitem{ZhuZhangYang2010} Zhu, S., Zhang, J., Yang, H.: Limiting profile of the blow-up solutions for the fourth-order nonlinear Schr\"odinger equation. Dyn. Partial Differ. Equ. 7(2), 187-205 (2010)

\bibitem{ZhuZhangYang2011} Zhu, S., Zhang, J., Yang, H.: Biharmonic nonlinear Schr\"odinger equation and the profile decomposition. Nonlinear Anal. 74(17), 6244-6255 (2011)

\bibitem{ZhuJEE2017} Zhu, S.: Existence of stable standing waves for the fractional Schr\"odinger equations with combined nonlinearities. J. Evol. Equ. 17(3), 1003-1021 (2017)

 \end{list}
\end{document}